\documentclass[11pt,leqno]{amsart}

\usepackage[utf8x]{inputenc}
\usepackage[T1]{fontenc}
\usepackage[english]{babel}

\usepackage{enumitem}

\usepackage[dvipsnames]{xcolor}
\usepackage{amsthm}
\usepackage[pagebackref, colorlinks=true, linkcolor=MidnightBlue, citecolor=OliveGreen]{hyperref}
\usepackage[capitalise,nameinlink]{cleveref}
\usepackage{amsmath,amssymb}
\usepackage{mathdots}
\usepackage{mathtools}
\usepackage{bbm}
\usepackage{thmtools}
\usepackage{thm-restate}

\usepackage{tikz}
\usetikzlibrary{automata, positioning}

\newtheorem{theorem}{Theorem}[section]
\newtheorem{lemma}[theorem]{Lemma}
\newtheorem{proposition}[theorem]{Proposition}
\newtheorem{corollary}[theorem]{Corollary}
\newtheorem{question}[theorem]{Question}

\theoremstyle{definition}
\newtheorem{definition}[theorem]{Definition}
\newtheorem{example}[theorem]{Example}

\AtBeginEnvironment{example}{%
  \pushQED{\qed}%
}
\AtEndEnvironment{example}{\popQED\endexample}

\DeclareMathOperator{\Z}{\mathbb{Z}}
\DeclareMathOperator{\N}{\mathbb{N}}

\DeclareMathOperator{\F}{\mathbb{F}}
\DeclareMathOperator{\Aut}{Aut}
\DeclareMathOperator{\Sym}{Sym}
\DeclareMathOperator{\Alt}{Alt}

\DeclareMathOperator{\St}{St}
\DeclareMathOperator{\st}{st}
\DeclareMathOperator{\Rist}{Rist}
\DeclareMathOperator{\rist}{rist}

\DeclareMathOperator{\id}{id}
\DeclareMathOperator{\syl}{syl}

\DeclareMathOperator{\reach}{reach}
\DeclareMathOperator{\core}{core}
\DeclareMathOperator{\ord}{ord}
\DeclareMathOperator{\spl}{fork}
\DeclareMathOperator{\Comp}{Comp}

\newcommand{\Dab}{\delta}

\renewcommand*{\backref}[1]{}
\renewcommand*{\backrefalt}[4]{%
  \ifcase #1 %
    No citations.
  \or
    (cited on page~#2).%
  \else
    (cited on pages~#2).%
  \fi%
}

\makeatletter
\@namedef{subjclassname@2020}{\textup{2020} Mathematics Subject Classification}
\makeatother

\title{Maximal subgroups in torsion branch groups}

\author[M. E. Garciarena]{Mikel Eguzki Garciarena}
\address{Mikel Eguzki Garciarena: Department of Mathematics, University of Salerno, 84084 Fisciano (SA), Italy; Department of Mathematics, University of the Basque Country UPV/EHU, 48080 Bilbao, Spain}
\email{mgarciarenaperez@unisa.it}

\author[J. M. Petschick]{J. Moritz Petschick}
\address{Jan Moritz Petschick: Fakult\"at für Mathematik, Universität Bielefeld, D-33501 Bielefeld, Germany}
\email{jpetschick@math.uni-bielefeld.de}

\thanks{}
\date{\today}

\begin{document}
\begin{abstract}
	We study the maximal subgroups of branch groups and obtain a criterion that ensures that certain spinal groups are contained in the class $\mathcal{MF}$ of groups with all maximal subgroups of finite index. This allows us to construct branch groups within $\mathcal{MF}$ exhibiting novel properties, for example groups that possess non-normal maximal subgroups. Furthermore, we give new concrete examples of branch groups outside $\mathcal{MF}$, with explicitly given maximal subgroups of infinite index. Most prominently, we construct a periodic branch group outside~$\mathcal{MF}$.
\end{abstract}

\maketitle

\section{Introduction} 
\label{sec:introduction}

\subsection*{Overview} 
\label{sub:overview}

{Branch groups are certain groups} of graph automorphisms of rooted trees sharing characteristics with the full automorphism group {of such a tree. They naturally occur} in the study of infinite groups with no proper infinite quotient. Moreover, they are known to be a rich source of examples and counter-examples in group theory, e.g.\ of groups of intermediate word growth \cite{Gri83}, of finitely generated infinite periodic groups \cite{Gri80, Pet23}, or of finitely generated non-nilpotent Engel groups \cite{Pet24}, {to only name a few instances}.

{The study of maximal subgroups of a group is a classical subject that is intimately connected with the study of the group's primitive actions. We are interested in groups that are highly non-primitive, namely infinite groups that only permit primitive actions on finite sets. Such groups are precisely those infinite groups whose maximal subgroups are all of finite index; we shall denote the class of such groups by $\mathcal{MF}$. Due to a classical theorem of Margulis and So\u ifer, the finitely generated linear groups in $\mathcal{MF}$ are precisely the virtually soluble ones. As a consequence, finitely generated groups permit a faithful primitive action if and only if they are of `almost simple type'; see~\cite{GGS22} for a detailed survey on groups with this property. This statement does not generalise to the case of general groups, as exemplified by Grigorchuk's group, a finitely generated branch (and hence non-linear) group that is easily seen to be of almost simple type, but is a member of $\mathcal{MF}$ by a result of Pervova \cite{Per05}. In this paper, we provide a machinery to produce many more such examples with various novel properties.

Independent from the above, Pervova's inquiry into the maximal subgroups of the Grigorchuk group and of the periodic GGS-groups --- another class of finitely generated branch groups --- was motivated by a conjecture of Kaplansky~\cite{Kap70} on the group algebra of a finitely generated group over a field of positive characteristic, proposing that if the Jacobson radical and the augmentation ideal of the group algebra coincide, the group in question is a finite $p$-group. Passman \cite{Pas98} reduced the conjecture to the case of (not necessarily finite) $p$-groups with every maximal subgroup normal of index $p$, infinite $p$-groups within $\mathcal{MF}$ such as those provided by Pervova are natural candidates for counter-examples. However, all branch $p$-groups that are known to be contained in $\mathcal{MF}$ do not meet the prerequisite of Kaplansky's conjecture.

There is a second conjecture, stated by Passman \cite[Problem 18.81]{KM20}, that is relevant to our subject, which asks if there exists an infinite residually finite $p$-group with a maximal, but not normal, subgroup. Since such a subgroup must necessarily be of infinite index, here infinite $p$-groups \emph{outside} $\mathcal{MF}$ are natural counter-examples.

As branch groups are a natural source of many infinite $p$-groups, it is thus of interest to determine when they fall into the class $\mathcal{MF}$ and when they do not; making them possible counter-examples to the conjectures of Kaplansky and Passman, respectively. Indeed, following Pervova's initial article, a long chain of research has dealt with the maximal subgroups of branch groups. Following a question formulated by Bartholdi, Grigorchuk and \v Suni\'k, the existence of finitely generated branch groups outside $\mathcal{MF}$ was non-constructively (i.e.\ without providing a concrete maximal subgroup of infinite index) confirmed by Bondarenko in \cite{Bon10}. Later, Francoeur and Garrido \cite{FG18} gave the first concrete example of a pair of branch group and a maximal subgroup of infinite index within, and even classified all maximal subgroups, showing that there are countably many of them, amongst them a finite number of finite index.} In the present work, we add to this list by providing a construction that allows to make Bondarenko's result constructive in certain cases.

On the other hand, Pervova's results have been generalised to larger families of (weakly) branch groups: to prime non-periodic \cite{FT22} and prime multi-GGS \cite{AKT16}, to the Basilica group \cite{Fra20} and relatives \cite{RT24}, to prime multi-EGS groups \cite{KT18}, and to GGS groups acting on trees of growing valency \cite{ST24}. {Furthermore, Bou-Rabee, Leemann and Nagnibeda \cite{BLN16} went on towards studying weakly maximal subgroups, i.e.\ infinite index subgroups not properly contained in any infinite index subgroup.

The goals of this paper are twofold. On one hand, we provide a machinery that allows to prove that certain members of a wide-ranging class of groups --- containing many branch groups --- are members of $\mathcal{MF}$, thereby far generalising previous results. This allows for the easy construction of many previously unknown groups without maximal subgroups of infinite index. On the other hand, we provide a novel approach to construct pairs of branch groups and their corresponding maximal subgroups of infinite index explicitly.  We use both methods to provide examples of new combinations of properties.

We consider the class of constant spinal groups, which are far generalisations of the class of GGS groups considered by Pervova and others. These are groups of the following form. Let $R$ be some finite group and denote by $R^\ast$ the infinite rooted tree consisting of words with letters in the group $R$. The group $R$ embeds into the automorphism group of $R^\ast$ by `rooted automorphisms', i.e.\ automorphisms that only modify the first letter of a word --- an element of $R$ --- by right-multiplication. Furthermore, the group $R^{|R|-1}$ also embeds into the automorphism group by so called `directed automorphisms', i.e.\ automorphisms which act by rooted automorphisms on all subtrees $rR^\ast$ for $r$ non-trivial and recursively by themselves on the subtree $eR^\ast$. A group generated by the group of rooted automorphisms and a group $D$ of directed automorphisms satisfying some additional property is called constant spinal. For a precise definition, see \cref{sub:constant_spinal_groups}. In the case of cyclic groups $R$ and $D$, one recovers the GGS groups.


\subsection*{Results on containment in $\mathcal{MF}$} 
\label{sub:results_on_containment_in_mathcal_mf}

The (finite) input data defining a constant spinal group allows the construction of a finite directed graph with its vertices given by members of the factor commutator group of the rooted group, from whose connectivity properties we may read of a criterion ensuring membership in $\mathcal{MF}$, which is recorded in \cref{thm:main positive}. Because of the extensive nature of its assumptions, we refrain from stating it at this point; however, the proof of \cref{thm:main positive} is the main technical achievement of this paper.}

Using \cref{thm:main positive}, {it is straight-forward not only to recover but to improve previous results on GGS groups to include groups acting on primary rather than prime valency.} Note that while the class of GGS groups is much more restricted than general constant spinal groups, it has received a great deal of attention, cf.\ \cite{FZ13,GGN24} for some in-depth structural results.

\begin{restatable}{thm}{primarymultiGGS}\label{thm:primary multi-GGS}
	Let $G = \langle R \cup D \rangle$ be a periodic primary multi-GGS group acting on the rooted tree $\underline{R}^\ast$. Assume that there exists a directed element $d \in D$ and a vertex $x \in \underline{R}$ corresponding to a generator of $R$, and such that $d|_{x}$ also generates $R$. Then $G$ is contained in $\mathcal{MF}$.
\end{restatable}

Here we distinguish between the finite group $R$ and its underlying set $\underline{R}$ from which the rooted tree is build to avoid ambiguities regarding the group operation and the concatenation of words over the alphabet $\underline{R}$. The notation $d|_{x}$ refers to the automorphism of the rooted subtree $x\underline{R}^\ast$ induced by $d$.

{Continuing, we use the fact that constant spinal groups permit arbitrary rooted groups to answer a question of Garrido and Francoeur on the existence of branch groups contained in $\mathcal{MF}$ but not in the class $\mathcal{MN}$ of groups with all maximal subgroups normal affirmatively, see \cref{thm:non mn}. To achieve this, we define a branch group contained in $\mathcal{MF}$ with a non-nilpotent rooted group. Since the rooted group is always (isomorphic) to a quotient of a constant spinal group, there exists a non-normal maximal subgroup. In contrast to our construction, all previous examples of branch groups inside $\mathcal{MF}$ did exclusively act cyclically on the first layer vertices.} Moreover, applying our method to a more sophisticated group, we prove the following somewhat stronger result; {here, a group is called just-insoluble if is insoluble, but every proper quotient is soluble.}

\begin{restatable}{thm}{mfnotjustinsol}\label{thm:mf not just-insol}
	There exists a finitely generated branch group within $\mathcal{MF}$ that is not just-insoluble.
\end{restatable}

{Note that both theorems are proven by explicitly constructing groups with the desired properties, using methods that may be adapted to achieve a wide variety of branch groups contained in $\mathcal{MF}$ exhibiting diverse features.}


\subsection*{Results on non-containment in $\mathcal{MF}$} 
\label{sub:results_on_non_containment_in_mathcal_mf}

{Turning} our attention to branch groups outside of $\mathcal{MF}$, {we make the following definition:} A special subgroup of a constant spinal group $G = \langle R \cup D \rangle$ is a subgroup that is conjugate to a subgroup of the form $\langle R \cup T \rangle$ for some subgroup $T \leq D$. We consider groups {of a type} that {forbids inclusion in} $\mathcal{MF}$ by the results of Bondarenko and construct, using special subgroups, concrete examples of maximal subgroups of infinite index. Here, the key tool is the following theorem.

\begin{restatable}{thm}{maximalsubgroupsinlayeredconstantspinalsubgroups}\label{thm:maximal subgroups in layered constant spinal subgroups}
	Let $G = \langle R \cup D \rangle$ be a constant spinal group, and let $H$ be a special prodense subgroup of $G$ conjugate to $\langle R \cup T \rangle$, which is branch over its commutator subgroup $H'$. The subgroup $H$ is a maximal subgroup of $G$ if and only if $T$ is a maximal subgroup of $D$.	
\end{restatable}

{To achieve this theorem, we establish some results on the structure of constant spinal groups that are interesting in their own right, cf.\ \cref{lem:labels of finitary elements} and \cref{lem:layer-climbing lemma}.}

Using the theorem above, we can construct many concrete examples of maximal subgroups of infinite index in branch groups. In particular, we achieve the following.

\begin{restatable}{thm}{torsionbranchgroupoutsidemf}\label{thm:torsion branch group outside mf}
	There exists a periodic branch group outside $\mathcal{MF}$.
\end{restatable}

Note that this is of special interest since \v Suni\'k groups --- which form the class of branch groups investigated by Francoeur and Garrido --- are known to be outside $\mathcal{MF}$ if they are non-periodic, and are conjectured to be included in $\mathcal{MF}$ if they are periodic. Our result shows that the existence of elements of infinite order is not a requirement for containment in $\mathcal{MF}$.

As another application, we can generalise Bondarenko's result on the existence of branch groups outside $\mathcal{MF}$ in the following way.

\begin{restatable}{thm}{arbitrarylongchainsofmaximalsubgroups}\label{thm:arbitrary long chains of maximal subgroups}
	Let $n \in \N$. There exists a branch group $G$ containing a chain of subgroups
	\[
		G_0 \leq G_1 \leq G_2 \leq \dots \leq G_n = G
	\]
	such that $G_{i-1}$ is a branch group that is a maximal subgroup of infinite index in $G_{i}$ for $i \in [n]$.
\end{restatable}

All previously known examples of branch groups outside $\mathcal{MF}$ possess the congruence subgroup property, i.e., every normal subgroup of finite index contains a (pointwise) stabiliser of a set of vertices of fixed distance~$n$ to the root of the tree. On the other hand, there are branch groups contained in $\mathcal{MF}$ with and without the congruence subgroup property \cite{KT18,Per07}. Our final result completes the picture.

\begin{restatable}{thm}{branchgroupoutsidemfwithoutcsp}\label{thm:branch group outside mf without csp}
	There exists a branch group outside $\mathcal{MF}$ without the congruence subgroup property.
\end{restatable}

{Along the way to prove the last result we give a characterisation of layered groups (i.e.\ groups that are branch over themselves) with the congruence subgroup property and use it to find a large family of branch groups without this property.}


\subsection*{Organisation} The paper follows the following structure: In \cref{sec:rooted_trees_and_constant_spinal_groups} we give some preliminaries, {mostly on groups acting on rooted trees. We define the class of constant spinal groups and give a description of their structure, building on previous results}. \cref{sec:constant_spinal_groups_with_all_maximal_subgroups_of_finite_index} contains {the development of our machinery to establish containment in $\mathcal{MF}$ and the} proof of \cref{thm:main positive}, which is then applied to prove \cref{thm:non mn}, \cref{thm:mf not just-insol} and \cref{thm:primary multi-GGS} in \cref{sub:applications_of_cref_thm_main_positive}. Finally, in \cref{sec:maximal_subgroups_of_infinite_index_in_branch_groups} we {further develop the structural theory on constant spinal groups and, building on this,} prove \cref{thm:maximal subgroups in layered constant spinal subgroups}, {which is afterwards applied to prove} \cref{thm:torsion branch group outside mf}, \cref{thm:arbitrary long chains of maximal subgroups} and \cref{thm:branch group outside mf without csp}, {along with some further examples of concrete subgroups of infinite index}.

\subsection*{Acknowledgements}
The first author is supported by the Spanish Government, grant PID2020-117281GB-I00, partly with FEDER funds and from the “National Group for Algebraic and Geometric Structures, and their Applications” (GNSAGA - INdAM). {The first author would also like to thank the Department of Mathematics at the Bielefeld University for its excellent hospitality and the support by the CRC-TRR~358 in Bielefeld.} The second author is supported by the Deutsche Forschungsgemeinschaft (DFG, German Research Foundation) – Project-ID 491392403 – TRR~358.


\section{Rooted trees and constant spinal groups} 
\label{sec:rooted_trees_and_constant_spinal_groups}

\subsection{Notation} 
\label{sub:notation}

Given integers $a, b \in \Z$, the intervals $[a, b], (a, b), (a,b]$ and $[a, b)$ denote the set of integers contained in the corresponding real intervals. Furthermore, we set $[n] = [1, n]$. For use in set valued maps, singleton sets $\{x\}$ are written just as their unique element $x$ without brackets. Given two finite lists $\epsilon \in X^n, \delta \in X^m$ from the same set $X$, the list $\epsilon \circ \delta$ denotes the list $\epsilon$ with $\delta$ appended. For a group~$G$, the natural epimorphism to the commutator factor group is denoted $\pi_G^\mathrm{ab}\colon G \to G/G'$.

\subsection{Graphs} A graph is a pair $(V, E)$ of a set of vertices and a set of edges $E \subseteq \binom V 2$, i.e.\ without multiple edges or loops. A graph automorphism is a map $g \colon V \to V$ such that $(g(v), g(u)) \in E$ for all $(v, u) \in E$.

A directed graph is a triple $(V, E, r)$ of a set $V$ of vertices, a set $E$ of edges, and a map $r\colon E \to V\times V$ assigning to every edge a tuple consisting of a source and target vertex. In particular, edges have an orientation and we allow loops and multiple edges. A directed graph is called simple if it contains no multiple edges, i.e.\ for every pair $v, w \in V$ of vertices there exists at most one edge $e \in E$ with $r(e) = (v,w)$. Note that simple directed graphs are still permitted to contain loops.

An (edge-)labelling of a directed graph $(V, E, r)$ is a map $\lambda\colon E \to L$ from $E$ to a set of labels $L$. A directed graph together with an labelling is called an (edge\nobreakdash-)labelled directed graph. Such a graph is called \emph{simply labelled} if the directed subgraphs $(V, \lambda^{-1}(l), r)$ are simple for every $l \in L$.

A path in a directed graph $(V, E, r)$ is a finite non-empty list of edges $(e_1, \dots, e_l) \in E^l$ such that the target of $e_i$ is equal to the source of $e_{i+1}$ for all $i \in [l-1]$. The length of a path is the length of the list. Given a simply labelled directed graph, paths are encoded by the sequence of the labels of their edges, with the starting point implicit.


\subsection{Rooted trees and their automorphisms} 
\label{sub:rooted_trees_and_their_automorphisms}

Let $X$ be a finite non-empty set which we call the \emph{alphabet}. We denote by $X^\ast$ the \emph{($X$-regular) rooted tree}, i.e.\ the Cayley graph of the free monoid on the set $X$. Every vertex is represented by a finite word $x_1 \dots x_n$ with $x_i \in X$ for $i \in [n]$ and $n \in \N_0$, which is said to be of \emph{length}~$n$, and there is an edge between $u$ and $v$ if and only if there exists $x \in X$ such that $ux = v$.  In this case the vertex~$u$ is called the \emph{predecessor} of the vertex~$v$. The unique word of length~$0$ (representing the neutral element) is denoted $\varnothing$ and is called the root. It is the only vertex of valency $|X|$, while every other vertex has valency $|X| + 1$. The~\emph{$n$\textsuperscript{th} layer} of the tree is the set $X^n$ of all words of length $n$, or equivalently, of all vertices of (geodesic) distance $n$ to $\varnothing$. We denote by $\Aut(X^\ast)$ the group of all graph automorphisms of $X^\ast$. We reserve the symbol $\id$ for the identity automorphism of~$\Aut(X^\ast)$. A \emph{ray} $u = (u_n)_{n \in \N_0}$ in $X^\ast$ is a sequence of vertices such that  $u_n$ is the predecessor of $u_{n+1}$ for every $n \in \N_0$. The set of rays is denoted by $\partial X^\ast$.

For $n \in \N$, we denote by $\St(n)$ the \emph{$n$\textsuperscript{th} layer} stabiliser, that is the pointwise stabiliser of the set~$X^n$. Similarly, for a vertex $v$ the point stabiliser is denoted $\st(v)$. For $g \in \Aut(X^\ast)$ and $v \in X^\ast$, we denote by $g|_v$ the \emph{section of $g$ at $v$}, that is the unique automorphism of $X^\ast$ such that we find $g(vu) = g(v) g|_v(u)$ for every vertex $u$. For a set $S \subseteq X^\ast$, we write $g|_S$ for the set $\{g|_s \mid s \in S\}$. We will make constant use of the following identities for sections under composition and inversion,
\begin{align*}
	(fg)|_u = f|_{g(u)} g|_u \quad \text{and} \quad f^{-1}|_u = (f|_{f^{-1}(u)})^{-1}.
\end{align*}
In particular, if $f \in \st(u)$ for some $u \in X^{\ast}$ we derive the following useful rule for conjugation,
\begin{align*}
	(f^g)|_u = (f|_{g(u)})^{g|_u}.
\end{align*}
Similarly, we denote by $g|^v$ the \emph{label of $g$ at $v$}, that is the unique permutation $g|^v \in \Sym(X)$ such that $g(vx) = g(v)g|^v(x)$ for every $x \in X$. The map $X^\ast \to \Sym(X)$ assigning to every vertex the corresponding label of $g$ is called the \emph{portrait} of $g$. Any automorphism is determined by its portrait. Let $R \leq \Sym(X)$ be a subgroup. We denote by $\Gamma(R)$ the subgroup of $\Aut(X^\ast)$ consisting on all automorphisms whose labels are elements of $R$.

The group $\Aut(X^\ast)$ naturally carries the structure of a wreath product $\Aut(X^\ast) \wr \Sym(X)$, where the base group is identified with $\St(1) \cong \Aut(X^\ast)^X$. The projections to the components of $\St(1)$ are given by the section maps $g \mapsto g|_x$ for $x \in X$. Write~$\psi \colon \Aut(X^\ast)^X \to \Aut(X^\ast)$ for the embedding of the base group into $\Aut(X^\ast)$; we will often suppress that symbol and think of $\Aut(X^\ast)^X$ as a subgroup of $\Aut(X^\ast)$.

We say that an element $g$ of $\Aut(X^\ast)$ is \emph{rooted} if $g|_x = \id$ for all $x \in X$,  i.e.\ if it is an element of the top group in the wreath product given above. Given a permutation $\sigma \in \Sym(X)$, by abuse of notation we simply write $\sigma$ for the rooted automorphism which acts as $\sigma$ on the vertices of the first layer.

We introduce the following shorthand notation. Let $n \in \N$, let $X_1, \dots, X_n$ be a collection of disjoint subsets of $X$ and let $g_1, \dots, g_n \in \Aut(X^\ast)$ be a collection of automorphisms. Then the symbol
\[
	(X_1: g_1, \, \dots \, X_n: g_n)
\]
stands for the unique automorphism $g \in \St(1)$ such that $g|_x = g_i$ for all $x \in X_i$ and~$i \in [n]$ and $g|_y = \id$ for all $y \in X \smallsetminus (\bigcup_{i = 1}^n X_i)$. If some of the sets $X_i$ are singletons~$\{x_i\}$, we simply write $x_i$.

By the last two paragraphs, any element $g \in \Aut(X^\ast)$ may be written as a product of the form
\begin{align*}
	g = (x: g|_x)_{x \in X} \sigma
\end{align*}
for some $\sigma \in \Sym(X)$.

For a subset $S \subseteq \Aut(X^\ast)$ and a vertex $u \in X^\ast$, we denote by $S_u$ the \emph{projection of $S$ at $u$}, that is the set 
\[
	S_u = \st_S(u)|_u.
\]


\subsection{Finite groups as alphabets} 
\label{sub:groups_as_alphabets}

Let $R$ be a finite group and write $\underline{R}$ for its underlying set. Using its Cayley embedding, we identify $R$ with a subgroup of $\Sym(\underline{R})$ acting by $r \underline{s} = \underline{rs}$. In this paper, we (mostly) consider subgroups of $\Gamma(R)$ acting on the $\underline{R}$-regular rooted tree. As a consequence, $r \in R$ is seen as a rooted automorphism, while $\underline{r}$ is a first layer vertex. In particular, the (notationally omitted) operation on underlined letters $\underline{r}$ and $\underline{s} \in \underline{R}$ is the product in the free monoid on $X = \underline{R}$, while the products $rs = t$ and $r\underline{s} = \underline{t}$ are computed by the group multiplication of $R$.

To reduce the amount of underlined letters, we write still write $X$ for $\underline{R}$ and denote first layer vertices by the letters $x$ or $y$, expect when it is useful to do otherwise.

Note that the restriction to trees built by finite groups is by no means an exotic restriction, as this covers e.g.\ the Sylow pro-$p$ subgroup given by $\Gamma(\langle \sigma \rangle)$ for any $p$-cycle $\sigma \in \Sym(X)$ for a set $X$ of cardinality $p$. Many (if not most) well-studied branch groups --- as for example the Grigorchuk or the Gupta--Sidki $p$-groups --- are subgroups of such a Sylow pro-$p$ subgroup.

The set $X = \underline{R}$ has a natural distinguished point, namely the trivial element $\underline{e}$ of the group~$R$. An element $g \in \Aut(X^\ast)$ is called \emph{directed} if $g|^v \neq e$ only if $v = \underline{e}^nx$ for some $n \in \N$ and $x \in X$. (Recall that $\underline{e}^n$ is the word $\underline{e} \dots \underline{e}$ of length $n$.)


\subsection{Groups of automorphisms} 
\label{sub:groups_of_automorphisms}

Let $G$ be a subgroup of $\Aut(X^\ast)$. We put $\st_G(v) = \st(v) \cap G$ and $\St_G(n) = \St(n)\cap G$, respectively, for all vertices~$v \in X^\ast$ and all $n \in \N$. We say that $G$ is \emph{self-similar}, if for every $g \in G$ and every $v \in X^\ast$, the section $g|_v\in G$ is contained in $G$. A self-similar group $G$ is called \emph{fractal} if $\st_G(v)|_v = G$ for all vertices $v \in X^\ast$. We say that {a subgroup} $G{\leq \Aut(X^\ast)}$ is \emph{spherically transitive} if it acts transitively each set $X^n$ for $n \in \N$.

Let $G$ and $H \leq \Aut(X^\ast)$ be two subgroups. We say that $H$ is \emph{geometrically contained in $G$} if $\psi(H^X) \leq G$, i.e.\ if $(x: h_x) \in G$ for every collection of elements $h_x \in H$.

Given a group $G \leq \Aut(X^\ast)$ and a vertex $v \in X^\ast$, the \emph{rigid vertex stabiliser of $v$} is the subgroup $\rist_G(v) \leq \st_G(v)$ consisting of all automorphisms possessing non-trivial labels only on the subtree $vX^\ast$ of vertices with prefix $v$. For $n \in \N$, we denote by $\Rist_G(n)$ the \emph{$n$\textsuperscript{th} rigid layer stabiliser}, which is defined as the product of all rigid vertex stabilisers of vertices of {layer}~$n$.

A spherically transitive group $G$ is called \emph{branch} if $\Rist_G(n)$ is a finite index subgroup for every $n \in \N$. {Let $K$ be a finite index subgroup of $G$.} We say that $G$ is \emph{regular branch {over $K$}} if $K$ is geometrically contained in $K \leq G$ as a finite index subgroup. If $G$ is self-similar, the geometric containment is always of finite index. If a self-similar group $G$ is regular branch over itself, we say that $G$ is \emph{layered}.

The following lemma is very useful to establish that some group is {regular branch}.

\begin{lemma}\label{lem:passage to normal closure}\cite[Proposition~2.18]{FZ13}
	Let $G$ be a spherically transitive {and} fractal group. Let $g \in G$ and $x \in X$. If the element $(x: g)$ is contained in $G$, then $\langle g \rangle^G$ is geometrically contained in $\langle (x: g) \rangle^G$.
\end{lemma}

{A subgroup $G \leq \Aut(X^\ast)$} is said to have the \emph{congruence subgroup property} if every { finite index} normal subgroup of $G$ contains some {layer} stabiliser{, i.e.\ if the completion $\overline{G}$ of $G$ with respect to the \emph{congruence topology} --- defined by taking the collection of {layer} stabilisers as a neighbourhood basis for the identity --- is isomorphic to the profinite completion of $G$.} A subgroup $H$ of $G$ is said to be \emph{congruence-dense} if {it is dense in the congruence topology, i.e.\ if $H\St_G(n) = G$ for every $n \in \N$.}


\subsection{Maximal subgroups of branch groups} 
\label{sub:maximal_subgroups_of_branch_groups}

{Recall that $\mathcal{MF}$ denotes the class of groups with all maximal subgroups of finite index. Let $G$ be a group. A subgroup $H \leq G$ is called \emph{prodense} if $HN = G$ holds for all non-trivial normal subgroups $N \trianglelefteq G$, i.e.\ if it is dense in the \emph{normal topology}, which is the topology of $G$ with the set of all non-trivial normal subgroups as a identity basis. The subgroup $H$ is called \emph{profinitely dense} if $HN = G$ holds for all normal subgroups $N \trianglelefteq_\mathrm{f} G$ of finite index, i.e.\ if it is dense in the profinite topology.}

For finitely generated groups, the existence of maximal subgroups of infinite index is tied to the existence of proper prodense subgroups{, as exhibited by the following result of Francoeur.}

\begin{lemma}\cite[Proposition~2.21]{Fra20}\label{lem:mf no proper prodense}
	Let $G$ be a finitely generated {infinite group such that every proper quotient of $G$ is a member of $\mathcal{MF}$.} Then $G$ is a member of $\mathcal{MF}$ if and only if $G$ contains no proper prodense subgroup.
\end{lemma}
The restriction to finitely generated groups is necessary, otherwise we may be in the situation that there exist proper subgroups that are not contained in any maximal subgroup at all, cf.\ \cite{Neu37}. Furthermore, if a proper quotient of $G$ does possess a maximal subgroup of infinite index, then so does $G$.

We now look at branch groups. By the lemma above, we may focus our attention on proper prodense subgroups. This has been the primary approach for studying maximal subgroups in branch groups, starting with the work of Pervova. {As a consequence,} we have some understanding of the structure of proper prodense subgroups in branch groups. In particular, we shall rely on a theorem of Francoeur.

\begin{theorem}\cite{Fra20}\label{thm:projections of proper pd}
	Let $G$ be a branch group acting on a rooted tree $X^\ast$, let $H \leq G$ be a prodense subgroup and let $u \in X^\ast$ be a vertex. Then $H_u$ is a prodense subgroup of $G_u$, and $H$ is a proper subgroup of $G$ if and only if $H_u$ is a proper subgroup of $G_u$.
\end{theorem}

See also \cite[Proposition~6.3]{FG18} for a similar result that only requires profinitely dense subgroups. Furthermore, due to another result of Francoeur, {every} maximal subgroup of infinite index in {a} branch group {is} again a branch group; our results in \cref{sec:maximal_subgroups_of_infinite_index_in_branch_groups} can be seen as a very explicit way to construct such pairs. Note that the converse does not hold, as the case of multi-GGS groups (which contain branch GGS groups as non-maximal subgroups of infinite index) exemplifies.

It is worth noting that may branch groups are \emph{just-infinite}, i.e.\ are groups that permit only finite proper quotients. For just-infinite groups, the normal and the profinite topologies naturally coincide. Thus every prodense subgroup is automatically profinitely dense. For further interplay between these properties in the general case, cf.\ \cite{GGS22}.

The just-infinite branch groups are described by the following theorem of Grigorchuk.

\begin{theorem}\cite[Theorem~4]{Gri00a}\label{thm:justinfinite}
	A branch group $G$ is just-infinite if and only if for each $v \in X^\ast$, the commutator factor group of $\rist_v(G)$ is finite. In particular, finitely generated periodic self-similar branch groups are just-infinite.
\end{theorem}

Applying this theorem, we arrive at the following easy corollary.

\begin{lemma}\label{lem:layered_csp_justinfinite}
	Let $G$ be a layered constant spinal group. Then $G$ is just-infinite.
\end{lemma}

\begin{proof}
	Since $G$ is layered, its rigid vertex stabilisers are isomorphic to $G$ itself. Since $G$ is spinal, it is finitely generated by element of finite order, hence $G/G'$ is finite. By \cref{thm:justinfinite}, $G$ is just-infinite.
\end{proof}


\subsection{Constant spinal groups} 
\label{sub:constant_spinal_groups}

Since { the occurrence of the first examples of branch groups in the 1980's, the original methods of construction have }been far generalised{, giving rise to} large and diverse families of groups, with some of their members exhibiting the same extraordinary properties as their inspiration {and some not}. { Particularly rich classes of such groups are given by the \emph{spinal groups} introduced in \cite{BGS03}. The object of our study are the so-called \emph{constant spinal groups}, which form a proper subclass of the class of spinal groups. They were previously studied by the second author in \cite{Pet23}. Here, we will slightly deviate from the definition given in \cite{Pet23} by restricting to groups of automorphisms whose labels are all acting regularly, cf.\ \cref{sub:groups_as_alphabets}. Our proofs build on this restriction at various places; however, it does not seem to be unavoidable in principle. Thus the following question seems to be of interest:

\begin{question}
	Does there exist a finitely generated branch group~$G$ in the class $\mathcal{MF}$ such that for every $u \in X^\ast$ there exists $g \in G$ such that $g|^u$ is non-trivial but not fixed point-free?
\end{question}}

\begin{definition}[Constant spinal group]\label{def:csg}
	Let $R$ be a finite group, put $X = \underline{R}$ and let $R$ act by left-multiplication on $X$. Let~$D \leq \Gamma(R)$ be a group of directed automorphisms such that
	\begin{equation*}\label{eq:generation property}
		D|_{X\smallsetminus\{e\}} = \langle d|_{x} \mid d \in D, x \in X\smallsetminus\{e\} \rangle = R. \tag{Gen}
	\end{equation*}
	A group~$G$ generated by $R \cup D$ is called a \emph{constant spinal group}. The group~$R$ is called the \emph{rooted group} of~$G$, and the group~$D$ is called the \emph{directed group} of~$G$.
\end{definition}


It is worth to note that $D$ is isomorphic to a subdirect product of $R$ (embedded into $R^X$). Thus $G/\St_G(2)$ is the subgroup of the wreath product $R \wr_X R$ generated by the top group and this subdirect product.

The most intensely studied family of constant spinal group is formed by the so-called (multi-)GGS groups{, the acronym standing for Grigorchuk--Gupta--Sidki}. Here, we may rely on previous results on the containment in $\mathcal{MF}$, wherefore we shall use them as a test case for our more general set-up.

\begin{definition}[GGS group]
	Let $G = \langle R \cup D \rangle$ be a constant spinal group. If $R$ is a cyclic group, we call $G$ a \emph{multi-GGS group}. If furthermore $D$ is cyclic, we {simply} speak of a GGS group. If $|X|$ is a prime or a prime power, respectively, we speak of \emph{prime}, or \emph{primary}, (multi-)GGS groups, respectively.
\end{definition}

Usually, multi-GGS groups are defined more directly, compare \cite{AKT16, FT22, FZ13, Vov00}. {To arrive at the standard definition, identify $X$ with the set $\Z/m\Z$ with $0$ corresponding to $\underline{e}$, fix a generator $a$ of $R$ and a minimal set of generators~$\{ b_{i} \mid i \in I \}$ for~$D$. For every $i \in I$ and $j \in X\smallsetminus\{0\}$, we find $b_{i}|_x = a^{e_{i;j}}$ for some $e_{i;j} \in \{0, \dots, m-1\}$. The tuples
\[
	e_i = (e_{i;1}, \dots, e_{i;m-1}) \in (\Z/m\Z)^{m-1}
\]
for $i \in I$ are usually called the defining tuples of~$G$. The requirement \eqref{eq:generation property} can be expressed as
\[
	\gcd \{ e_{i;j} \mid i \in I, j \in X\smallsetminus\{0\} \} = 1.
\]}
We now point out some basic properties shared by all constant spinal groups that we shall use constantly {and without further mention}. A proof may be found in \cite{Pet23}.

\begin{lemma}[Properties of constant spinal groups]\label{lem:properties csg}
	Every constant spinal group is spherically transitive, self-similar and fractal.
\end{lemma}

We remark that it is unclear to what extend the class of constant spinal groups {consists of branch groups}. The authors are only aware of a few examples of non-branch constant spinal groups ({which are still `weakly' branch, meaning that all rigid {layer} stabilisers are non-trivial}); these are the GGS groups with constant defining tuple, cf. \cite{FGU17}, and some $2$-groups introduced in \cite{Sie08}. All other prime multi-GGS groups on trees and most primary GGS groups are known to be branch {by~\cite{FZ13} and \cite{DFG22}, respectively}, and many other examples of constant spinal groups are indeed branch groups, see for example \cite{GS83a, Pet23}. However, no general criterium seems to be available.

It is readily observed that the property \eqref{eq:generation property} may already be achieved with a proper subgroup $T$ of $D${; see \cref{eg:class_of_examples} for a large class of examples.} In that case, the subgroup $H = \langle R \cup T \rangle$ of $\langle R \cup D \rangle$ is also a constant spinal group. Conjugates of subgroups of this form we call \emph{special}, they will play a major role in \cref{sec:maximal_subgroups_of_infinite_index_in_branch_groups}.

Fix a constant spinal group $G = \langle R \cup D \rangle$. Apart from the `natural' generating set~$R \cup D$, another, non-minimal generating set of $G$ is given by $R \cup D^R$. We say a word~$\gamma$ over the alphabet $R \cup D^R$ is of \emph{syllable form} if it fulfils
\[
	\gamma = r_0 \prod_{i = 1}^n d_i^{r_i}
\]
for some $n \in \N$ called the \emph{length}, some $r_i \in R$ and some $d_i \in D$ with $r_i \neq r_{i+1}$ for $i \in [1, n)$. The minimal integer $n$ among all such products evaluating to a given element $g \in G$ is called the \emph{syllable length} of $g$ and is denoted $\syl(g)$; any corresponding expression is called a \emph{minimal syllable form of $g$}. Thus the syllable length is the (weighted) word length with respect to the generating set $R \cup D^R$ with rooted elements having weight~$0$ and elements of $D^R$ --- that we shall call \emph{syllables} --- having weight~$1$. Note that the corresponding distance function is not a metric, but only a pseudo-metric on $G$.

Let $x \in X$ and $g \in G$. If $\gamma = r_0 \prod_{i = 1}^n d_i^{r_i}$ is a syllable form evaluating to $g$, we may find a syllable form $\gamma|_x$ evaluating to $g|_x$ in the following way. First compute
\[
	g|_x = \gamma|_x = (r_0 \prod_{i = 1}^n d_i^{r_i})|_x = \prod_{i = 1}^n d_i|_{r_i.x}.
\]
Notice that $d_i|_{r_i.x} \in D$ if and only if $r_i.x = \underline{e}$, otherwise $d_i|_{r_i.x} \in R$. Write $J_x(\gamma)$ for the set $\{i \in [1, n] \mid \underline{r_i^{-1}} = x\}$. Returning to the product $\prod_{i = 1}^n d_i|_{r_i.x}$, we make the following transformations: collect all rooted elements to the left, i.e.\ conduct the following operation as often until it is no longer possible to perform,
\[
	d r \rightarrow r d^{r},
\]
and then combine syllables in the same conjugate of $D$, i.e.\ conduct
\[
	d_1^r d_2^r \rightarrow (d_1d_2)^r.
\]
The resulting product
\[
	\gamma|_x = s_0 \prod_{i \in J_x(\gamma)} d_i^{s_i}
\]
is again of syllable form. It is called the \emph{induced syllable form with respect to $x$}. More generally, we set $\gamma|_u = \gamma|_{x_1}|_{x_2} \dots |_{x_n}$ for $u = x_1 \dots x_n \in X^n$.

Note that if $\gamma$ is of length $n$, the set $J_x(\gamma)$ (for any $x$) has at most cardinality $\lceil n/2 \rceil$ by the condition that neighbouring syllables are contained in different conjugates of~$D$. Thus the length of $\gamma|_x$ is at most $\lceil n/2 \rceil$. Furthermore, since the sets $J_x(\gamma)$ form a partition of $X$, we find that the overall length of all first {layer} sections is at most the length of $\gamma$. We record these well-known facts (cf.\ e.g. \cite[Lemma~2.4]{Pet23}) and its immediate consequences for the syllable length of an element in a lemma.

\begin{lemma}[Contraction lemma]\label{lem:contraction lemma}
	Let $G = \langle R \cup D \rangle$ be a constant spinal group and let $x \in X$.
	\begin{enumerate}
		\item Let $\gamma$ be a syllable form of length $n$. Then $\gamma|_x$ is of length at most $\lceil n/2 \rceil$.
		\item Let $g \in G$. Then $\syl(g|_x) \leq \lceil \syl(g)/2 \rceil$.
		\item Let $g \in G$. Then $\sum_{x \in X} \syl(g|_x) \leq \syl(g)$.
	\end{enumerate}
\end{lemma}

The second property is usually called \emph{contraction}. Since it allows to use induction on the length of elements, it is fundamental for the theory of constant spinal groups. More generally, {any} self-similar group is called contracting if the set
\[
	\mathfrak{N}(G) = \bigcup_{g \in G} \bigcap_{n \in \N} \{ g|_{uv} \mid u \in X^n, v \in X^\ast\},
\]
{which is} called the \emph{nucleus of $G$}, is finite, cf.\ \cite{Nek05}. Using the contraction lemma, we find $\mathfrak{N}(G) \subseteq R \cup D^R$ for every constant spinal group~$G$. Looking at the sections of those elements, we {furthermore obtain}
\[
	\mathfrak{N}(G) = R \cup D.
\]
The following result is essential for {our} analysis of constant spinal groups {that} are contained in $\mathcal{MF}$, {since we are frequently compute in the factor commutator group of~$G$ rather than in~$G$ itself.} Both the statement and its proof are adapted from the more general setting of polyspinal groups, which will be part of a forthcoming article. For completeness, we give a proof.

\begin{theorem}\label{thm:factor commutator groups of constant spinal groups}
	Let $G = \langle R \cup D \rangle$ be a constant spinal group. Then
	\[
		G/G' \cong R/R' \times D/D'.
	\]
\end{theorem}

\begin{proof}
	Since $G$ is generated by $R$ and $D$, its factor commutator subgroup is a quotient of $R/R' \times D/D'$. Evidently there is a surjective homomorphism $\rho \colon G \to G/\St_G(1) \cong R \to R/R'$. We will now construct a homomorphism $\delta$ from $G$ onto $D/D'$ with $R \leq \ker(\delta)$ by
	\[
		\delta(g) = \lim_{n \to \infty} \prod_{u \in X^n} \mathbf{1}_D(g|_u)D',
	\]
	where $\mathbf{1}_D$ is the function
	\begin{align*}
	\mathbf{1}_D(g)=\begin{cases} g\qquad&\text{if $g\in D$;}\\
	\id\qquad&\text{otherwise.}
	\end{cases}
	\end{align*} Let us first argue why this function is well-defined. By the description of the nucleus above, for every $g$ there exists some $n \in \N$ such that all sections at $u \in X^n$ are either rooted or directed. The rooted sections do not matter due to the definition of $\mathbf{1}_D$. Every directed section has one directed section that is furthermore equal to itself. Thus, the sequence $(\prod_{u \in X^n} \mathbf{1}_D(g|_u)D')_{n \in \N}$ is eventually constant and the limit well-defined.
	
	{We now} prove that $\delta$ is a homomorphism. Let $g, h \in G$. There exists some $n \in \N$ such that both $g$ and $h$ have only directed and rooted sections at vertices of layer $n$. Thus
	\[
		\delta(g) = \prod_{u \in X^m} \mathbf{1}_D(g|_u)D' \quad\text{and}\quad \delta(h) = \prod_{u \in X^m} \mathbf{1}_D(h|_u)D'
	\]
	for all $m \geq n$. Consider the product $gh$. Since $(gh)|_u = g|_{h(u)}h|_u$, all sections at vertices $u \in X^n$ are either rooted, directed or the product of a rooted and a directed element. {Partition $X^n$ into the disjoint union $(X^n)_{\mathrm{root}} \sqcup (X^n)_{\mathrm{dir}} \sqcup (X^n)_{\mathrm{prod}}$, collecting the vertices $u$ such that $(gh)|_u$ is rooted, non-trivial and directed, and neither, respectively.}

	For the computation of $\delta(gh)$, we may ignore the rooted sections of~$gh$. Every directed section is the product of two directed sections of $g$ and $h$ (where one of them may be trivial), and there is precisely one descendant of~$u$ where~$gh$ has the same directed section, all other sections at descendants of~$u$ are rooted. Thus
	\begin{align*}
		\prod_{v \in (X^n)_{\mathrm{dir}}X} \mathbf{1}_D((gh)|_v)D' &= \prod_{v \in (X^n)_{\mathrm{dir}}X} \mathbf{1}_D(g|_{h(v)}) \mathbf{1}_D(h|_v)D' \\&=
		\prod_{v \in h((X^n)_{\mathrm{dir}}X)}\left( \mathbf{1}_D(g|_{v})D' \right) \prod_{v \in (X^n)_{\mathrm{dir}}X}\left(\mathbf{1}_D(h|_v)D' \right).
	\end{align*}
	Now consider the case that $u \in (X^n)_{\mathrm{prod}}$, i.e.\ that either $g|_{h(u)}$ is rooted and $h|_u$ directed or the other way around. In either way, there is again precisely one descendant~$v$ of~$u$ such that $(gh)|_v$ is directed, and is in fact equal to the directed element among $g|_{h(u)}$ and $h|_u$. Thus we may compute as above. Putting everything together, we find
	\begin{align*}
		\prod_{v \in X^{m}} \mathbf{1}_D((gh)|_v)D' &=
		\prod_{v \in h(X^{m})}\left( \mathbf{1}_D(g|_{v})D' \right) \prod_{v \in X^{m}}\left( \mathbf{1}_D(h|_v)D' \right)
	\end{align*}
	for all $m > n$. Thus $\delta$ is a homomorphism. Considering the image of $D$ under $\delta$, it is plain to see that $\delta$ is surjective. Furthermore, for every rooted element $r$ we see $\delta(r) = D'$, hence ${(\rho \times \delta)}(r d) = (rR', dD')$, yielding a surjective homomorphism of~$G$ onto $R/R' \times D/D'$.
\end{proof}

In \cref{sec:constant_spinal_groups_with_all_maximal_subgroups_of_finite_index}, we shall consider constant spinal groups {that are periodic, i.e.\ such that every element has finite order}. It is not { completely} understood which constant spinal groups are periodic; however, it is well-known how to construct a plethora of periodic examples. We will employ a condition ensuring periodicity established by the second author in \cite{Pet23}, where the topic is treated in greater generality. To state that condition, we need some further terminology.

Let $g \in \Aut(X^\ast)$ be some automorphism of $X^\ast$ and let $u \in X^\ast$ be a vertex. Evidently, the orbit of $u$ under $g$ is finite, thus there exists some minimal $n \in \N$ such that $g^n \in \st(u)$. The \emph{stabilised section of $g$ at $u$} is defined to be
\[
	g\|_u = g^n|_u.
\]
Observe that for any $g \in G$ and $u \in X^\ast$, the stabilised section $g\|_x$ is naturally contained in the projection~$G_u$. It turns out that control over the length of stabilised sections of the elements of a group allows to establish periodicity. Given a constant spinal group~$G$, we refer to the set
\[
	\mathfrak{SN}(G) = \bigcup_{g \in G} \bigcap_{n \in \N} \{ g\|_{uv} \mid u \in X^n, v \in X^\ast \},
\]
as the \emph{stabilised nucleus} of~$G$; this concept is a variation of the nucleus of a self-similar group as stated above. Using it, we are able to state the following result.

\begin{proposition}\label{prop:periodicity via stab nuc}
	Let $G = \langle R \cup D \rangle$ be a constant spinal group. If $\mathfrak{SN}(G)$ consists of elements of finite order, the group $G$ is periodic.
\end{proposition}

The proof of this result is found in \cite[Proposition 2.8]{Pet23}, using a slightly different formulation; conceptually, it is an evolution of the ideas of Grigorchuk, Gupta and Sidki used to establish periodicity of their eponymous groups.

While it may be feasible to compute the stabilised nucleus of a given constant spinal group, it is difficult to obtain a useful description in general. Clearly, the stabilised nucleus of a constant spinal group contains the directed group, however, there may well be more elements in $\mathfrak{SN}(G)$. One possibility is to analyse the function $\Sigma_T \colon \mathcal{P}(R) \to \mathcal{P}(R)$ defined in the following for a constant spinal group $G = \langle R \cup D \rangle$. Here, $T$ denotes a generating set of the directed group $D$. For any $S \subseteq R$, we {put}
\begin{align*}
	\Sigma_T(S) = \bigcup_{r \in S} \bigcup_{s \in r^R} \langle \prod_{i = 1}^{\ord(s)-1} d|_{s^i} \mid d \in T \rangle \cdot \langle d|_{s^i} \mid i \in [1, \ord(s)), d \in T \rangle'.
\end{align*}
Note that the order of the products generating the left subgroup in the product defining~$\Sigma_T(S)$ is irrelevant due to the fact that the right subgroup is the derived subgroup of {the group generated by their factors}. Now we are able to state the following result, which will be chiefly of use to establish \cref{thm:non mn} and \cref{thm:mf not just-insol}.

\begin{theorem}{\cite[Theorem~B]{Pet23}}\label{thm:condition for stab nuc}
	Let $G = \langle R \cup D \rangle$ be a constant spinal group. Let~$T$ be a generating set for~$D$. If for all $r \in R$ there exists some $n \in \N$ such that $\Sigma_T^n(r) = \{\id\}$, then $\mathfrak{SN}(G) = D$, and, in particular, $G$ is periodic.
\end{theorem}



\section{Constant spinal groups with all maximal subgroups of finite index} 
\label{sec:constant_spinal_groups_with_all_maximal_subgroups_of_finite_index}

\subsection{General strategy to establish membership in $\mathcal{MF}$} 
\label{sub:general_strategy_to_establish_membership_in_mathcal_mf}

{Here we} give an overview {of the} general strategy {used} to prove that a contracting branch group is contained in the class~$\mathcal{MF}$. This method was first {developed} by Pervova in \cite{Per05}, and subsequently refined by Klopsch and Thillaisundaram \cite{KT18} and Francoeur \cite{Fra20}. {In all these cases, the methods only applied to groups contained a Sylow pro-$p$ subgroup of the automorphism group of a regular tree on a set of prime order.} We adapt and expand it {to cover groups whose labels do not necessarily commute; indeed, they may generate any finite group.}

{Our procedure is the following.} Fix a constant spinal group $G = \langle R \cup D \rangle$. Let~$M$ be a prodense subgroup of $G$. If one can show that there is a projection $M_u$ of~$M$ that is equal to $G_u = G$, the subgroup is not proper by \cref{thm:projections of proper pd}. Thus we analyse the projections of $M$. Given $u \in X^\ast$ and $m \in M$, the stabilised section $m\|_u$ {is} contained in the projection~$M_u$. Using the contraction of the syllable length under taking sections, we aim to find elements of length at most $1${ --- i.e. generators --- }in a projection. We repeat this step multiple times:
\begin{enumerate}
	\item Find some directed element $d$ in a projection. Whenever $d \in M_u$, it is also contained in every projection of $M_u$. This allows us to `keep' $d$ whenever we further descent to projections to lower vertices in the following steps. The first step is carried out in \cref{sub:obtaining_a_directed_element}.
	\item Since $M_u$ is prodense, for any $r \in R$ there exists some $g \in G'$ such that $rg \in M_u$. We describe a way to associate to every such pair $(r, g)$ a new pair $(r', g')$ such that $r'g' \in M_{ue}$, basically by considering the section $d^{rg}|_e$. The passage $r \mapsto r'$ is thus governed by the structure of sections of $d$, which we encode in a directed graph. Under some assumption on that graph, we understand that we are able to make sure that the collection of possible $r'$ that may be obtained generated the rooted group $R$. The second step is carried out in \cref{sub:definition_of_prodsec_n_d}.
	\item The main technical step is to analyse the behaviour of the syllable length under the passage from $g$ to $g'$. Under certain circumstances, we prove that the length eventually decreases, hence that there exists a projection $M_v$ of $M$ such that $M_v$ contains $d$ and the attainable rooted elements, which generate $R$ as a subgroup of $M_v$. The third step is carried out in \cref{sub:delta_and_theta_maps}, \cref{sub:theta_maps} and \cref{sub:reduction_along_paths}.
	\item Finally, using that all rooted elements are contained in $M_v$, one is able to gather more directed elements $t \neq d$ { by yet again passing to sections}. If one can obtain a generating set, we conclude that there exists a non-proper projection of $M$, which was our goal. The last step is carried out in \cref{sub:obtaining_a_directed_element}.
\end{enumerate}
Finally, we put everything together {and} prove \cref{thm:main positive} in \cref{sub:proof_of_cref_thm_main_positive}.


\subsection{Obtaining a directed element} 
\label{sub:obtaining_a_directed_element}

As a first step, we describe conditions under which every prodense subgroup $M \leq G$ has a projection $M_u$ containing a non-trivial directed element. We follow the strategy described in \cite{KT18}, which heavily uses that the group in question is periodic. Note that in the {special} case of GGS~groups, analogous statements {to many of our results in this section} for non-periodic groups {were obtained} in \cite{FT22} using different methods.

\begin{lemma}\label{lem:additivity of sections measured in the abelianisation of D}
	Let $G = \langle R \cup D \rangle$ be a constant spinal group. For all $g \in G$, the following equality holds true,
	\[
		\prod_{t \in T} \Dab(g\|_{\underline{t}}) = \prod_{x \in X} \Dab(g|_x) = \Dab(g),
	\]
	where $T$ denotes a transversal of the subgroup of $R$ generated by {$ g|^\varnothing $}.
\end{lemma}

\begin{proof}
	Reviewing the definition of $\Dab$, the second equality becomes clear. For the first one, write $rh = g$ for some $h \in \St_G(1)$ and $m$ for the order of $r$. Then consider that
	\[
		g\|_x = (rh)^m|_x = h|_{r^{m-1}x} \dots h|_{rx} h|_x = g|_{r^{m-1}x} \dots g|_{rx} g|_x.
	\]
	Thus every stabilised section $g\|_x$ is the product of all sections at vertices contained in the coset $\langle r \rangle x$, and we may conclude the proof by considering the following,
	\[
		\prod_{t \in T} \Dab(g\|_{\underline{t}}) = \prod_{t \in T} \prod_{i = 0}^{m-1} \Dab(g|_{r^i\underline{t}}) = \prod_{x \in X} \Dab(g|_x).\qedhere
	\]
\end{proof}

\begin{proposition}\label{prop:prodense conts directed elm}
	Let $G = \langle R \cup D \rangle$ be a constant spinal group such that $\mathfrak{SN}(G) = D$ and let $M \leq G$ be a prodense subgroup. Let ${\eta} \colon D/D' \to \Z/p^n\Z$ be a surjective homomorphism, where $p$ denotes some prime number and $n \in \N$ is some integer. Put ${\varphi = \eta \circ \Dab}$. There exists $u \in X^\ast$ such that ${\varphi}(M_u \cap D) = \Z/p^n\Z$.
\end{proposition}

\begin{proof}
	Let $g = rh$ be such that $h \in \St(1)$ and ${\varphi}(g) \not\equiv_p 0$. Since ${\varphi}$ factors through $\Dab$, we find that also for a transversal $T$ of $\langle r \rangle$ within $R$,
	\[
		\sum_{t \in T} {\varphi}(g\|_{\underline{t}}) = {\varphi}(g) \not\equiv_p 0.
	\]
	Consequently, there exists some $t \in T$ such that ${\varphi}(g\|_{\underline{t}}) \not\equiv_p 0$.
	
	Now, using surjectivity, let $d \in D$ be some element such that ${\varphi}(d)$ is not a multiple of~$p$. Since $M$ is prodense and $G'$ a normal subgroup of finite index, we find $MG' = G$, hence there exists some $c \in G'$ and $g \in M$ such that $gc = d$. Since ${\varphi}$ factors through~$G/G'$, we find ${\varphi}(g) = {\varphi}(d)$.
	
	{Overall}, we see that there exists some $x \in X$ such that ${\varphi}(g\|_x)$ is not a multiple of~$p$, {while} at the same time, $g\|_x \in M_x$ and, {by} the definition of $\mathfrak{SN}(G)$, either $g\|_x \in \mathfrak{SN}(G)$ {or $\syl(g) > \syl(g\|_x)$. Thus we may repeatedly} replace $g$ {by} $g\|_x$ until $g\|_u\in \mathfrak{SN}(G)$ for some $u \in X^\ast$. {Since} $\mathfrak{SN}(G) = D$ {by assumption}, we find $g\|_u \in M_u \cap D$ for some $u \in X^\ast$ with ${\varphi}(g\|_u) \not\equiv_p 0$. Evidently, the image of $g\|_u$ generates~$\Z/p^n\Z$.
\end{proof}

At this point, it is convenient to establish a proposition used for the {final} step in the general plan {laid out in \cref{sub:general_strategy_to_establish_membership_in_mathcal_mf}}: to locate of more than one directed element in a prodense subgroup under the assumption of the presence of sufficiently many rooted elements. To do so, we introduce the following.

\begin{definition}\label{def:core of a subgroup}
	Let $G \leq \Aut(X^\ast)$ be a group of automorphisms. The \emph{{core} of $G$} is the set
	\[
		\core(G) = \bigcup_{n \in \N}\bigcap_{N > n} G_{\underline{e}^N}.
	\]
\end{definition}

{Of course, there is nothing inherently distinctive about the constant ray $\overline{\underline{e}} = \{\underline{e}^N \mid N \in \N\}$; it is easy to adapt the above to definition to measure projections along any ray $u$ in the tree~$X^\ast$. However, since we define constant spinal groups using elements that are directed along the ray $\overline{\underline{e}}$, this definition fits our precise needs.

We now state some immediate consequences from the definition.}

\begin{lemma}
	Let $G \leq \Aut(X^\ast)$ be a group of automorphisms.
	\begin{enumerate}
		\item Let $n \in \N$. Then $\core(G) = \core(G_{\underline{e}^n})$.
		\item Let $H \leq \Aut(X^\ast)$ be another group of automorphisms. Then $\core(G \cap H)\leq \core(G) \cap \core(H)$.
		\item Let $D$ be a group of directed elements. Then $\core(D) = D$.
	\end{enumerate}
\end{lemma}

\begin{proposition}\label{prop:all cosets of dirs}
	Let $G$ be a constant spinal group with $\mathfrak{SN}(G) = D$. Let $M \leq G$ be a prodense subgroup such that $R \leq \core(M)$. Then $\Dab(\core(M) \cap D) = D/D'$.
\end{proposition}

\begin{proof}
	Since $R$ is a finite subgroup {of $\core(M)$}, there exists some integer $n \in \N$ such that $R \leq M_{e^m}$ for all $m \geq n$. Since $\core(M) = \core(M_{e^n})$ for all $n \in \N$, and since projections of prodense subgroups {remain} prodense, we may without loss of generality replace $M$ by $M_{e^n}$ and assume that $R \leq M_{e^n}$ for all $n \in \N$.
	
	Let $\overline{d} \in D/D'$. Since $M$ is prodense, we have $G = G'M$. Thus there exists some $g \in M$ such that $\Dab(g) = \overline{d}$. Write $r = g|^{\varnothing} \in R$. Let $T$ be a transversal of $\langle r \rangle$ in $X$. By \cref{lem:additivity of sections measured in the abelianisation of D}, we find
	\(
		\prod_{t \in T} \Dab(g\|_{\underline{t}}) = \Dab(g)
	\),
	while $g\|_{\underline{t}} \in M_{\underline{t}} = (M^{\underline{t}})_{\underline{e}} = M_{\underline{e}}$, using that $t \in R \leq M$. Since $R$ is contained in every projection $M_{\underline{e}^n}$ for $n \in \N$, we may repeat this argument arbitrary often and obtain a subset $T_n$ of $X^n$ such that $g\|_{\underline{t}} \in M_{\underline{e}^n}$ for all $t \in T_n$ and
	\[
		\prod_{t \in T_n} \Dab(g\|_{\underline{t}}) = \Dab(g).
	\]
	{Choosing a sufficiently large integer~$n$, we may assume --- by the definition of the stabilised nucleus --- that} $g\|_{\underline{t}}\in \mathfrak{SN}(G)$ { for all $t \in T$. By assumption, this allows us to conclude that all $g\|_{\underline{t}}$ are directed elements.} Thus, $h = \prod_{t \in T_n} g\|_{\underline{t}} \in M_{\underline{e}^n} \cap D$ satisfies $\Dab(h) = \Dab(g) = \overline{d}$. Finally, $h|_{\underline{e}^m} = h \in M_{\underline{e}^{n+m}} \cap D$ for all $m \in \N$, i.e.\ $h \in \core(M_{\underline{e}^n} \cap D)$. It remains to see that
	\[
		\core(M_{\underline{e}^{n}} \cap D) \leq \core(M_{\underline{e}^{n}}) \cap \core(D) = \core(M) \cap D.\qedhere
	\]
\end{proof}


\subsection{Filling directed elements} 
\label{sub:definition_of_prodsec_n_d}

We now describe {under} which circumstances a passage as described in the second step of the general plan may be feasible at all, i.e.\ allowing for recursion.

Every element $d \in D$ defines {---} and {may} indeed {be recovered from} {---} the directed graph, {whose} set of vertices is $R$, and {which possesses} an edge from $r \in R$ to $s \in R$ if and only if $d|_{\underline{r}} = s$. All sections of $d$ (except the {necessarily trivial} section at~$\underline{e}$) may be read of this graph.

However, when passing from a product $rg$ with $r \in R$ and $g \in G'$ to some appropriate $\tilde{r}\tilde{g}$ of the same form (as is our plan, see above), our information on $\tilde{r}$ as an element of $R$ is only accurate up to the commutator subgroup~$R'$. To properly capture this passage, we define a certain variation of the natural graph described above.

First, we establish another convention for convenience. For any $r \in R$ write
\[
	m(r) = \ord_{R/R'}(rR')
\]
for the order of the image of~$r$ in the {factor commutator group} $R/R'$. 

Now define a directed edge-labeled graph $\mathcal{R}^{\mathrm{ab}}(d)$, whose set of vertices is the set of cosets $R/R'$. Let $rR' \in R/R'$ be a vertex. Then there exists an edge labeled {$l$}, for some $l \in [1, m(r))$, from $rR'$ to the vertex $sR'$, if and only if there exists some $t \in R'$ such that $d|_{\underline{(rt)^l}}\equiv_{R'} s$.

The situation when every vertex of $\mathcal{R}^{\mathrm{ab}}(d)$ has out-valency equal to $m(r)-1$, {i.e.\ when the directed graph is simply labelled}, is {of particular interest to us}. This is the case whenever the $R'$-coset of a section at some $x \in X$ does only depend on $xR'$, i.e.\ when $d$ satisfies the following condition.

\begin{definition}\label{def:compatible directed elements}
	Let $G = \langle R \cup D \rangle$ be a constant spinal group. A directed element $d \in D$ is called \emph{compatible} if for all $r \in R$ and all $s \in R'$
	\[
		d|_{\underline{rs}} \equiv_{R'} d|_{\underline{r}}.
	\]
	The subgroup of all compatible elements is denoted $\Comp(D)$.
\end{definition}

Note that when $d$ is compatible, every $D$-conjugate is also compatible. Furthermore, since the sections of any element $t \in D'$ are contained in $R'$ via $[d_1, d_2]|_x = [d_1|_x, d_2|_x]$, if some $d$ is compatible, the whole coset $dD'$ is compatible, whence $\Comp(D)$ is a normal subgroup in $D$ fulfilling $D' \leq \Comp(D)$.

Using the graph $\mathcal{R}^{\mathrm{ab}}(d)$, we shall give a definition of when a directed element $d$ is {`filling' (cf.\ \cref{def:filling})}, i.e.\ allows us to recover $R$ from its sections whenever $d$ is contained in a prodense subgroup. To do so, we need to establish some graph-theoretical terminology.

\begin{definition}
	Let $\mathcal{G} = (V, E)$ be a directed graph. Let $v, w \in V$ be vertices of that graph. We say $w$ is \emph{reachable} from $v$ if there exists a directed path starting at $v$ and ending at $w$. Let $W \subseteq V$ be a set of vertices. Then the $\emph{reach}$ of $W$ is the set
	\[
		\reach(W) = \{ v \in V \mid v \text{ is reachable from } w \text{ for some } w \in W\}.
	\]
	In case $W = \{w\}$ is a singleton set, we just write $\reach(w)$ for $\reach(W)$.
\end{definition}

\begin{definition}
	Let $G = \langle R \cup D \rangle$ be a constant spinal group. Let $d {\in \Comp(D)}$ be a compatible directed element. We consider the graph $\mathcal{R}^{\mathrm{ab}}(d)$ as defined above. The set $\spl(d)$ of \emph{forking points} consists of the set of {vertices} $rR' \in R/R'$ {satisfying} the following {list of} conditions:
	\begin{enumerate}
		\item there exists a directed path of length $\ell$, starting and ending in~$rR'$, whose first edge is labeled with $k \in [m(r)-1]$ such that $k$ is coprime to the order of $d$, and
		\item there exists a directed path of length $\ell'$, starting and ending in~$rR'$, whose first edge is labeled with $k' \in [m(r)-1]$ such that $k' \neq k$, and
		\item the lengths $\ell$ and $\ell'$ are coprime.
	\end{enumerate}
	The pair of integers $(k, k')$ is called the \emph{signature} of the forking point $rR'$.
\end{definition}

It is easily seen that $m(r) > 2$ for all forking points of some $d \in D$. Note that it is acceptable that the length $\ell$ of one of the required paths is $1$. By the next lemma, we are allowed to restrict our attention to forking points with signature $(1, l)$.

\begin{lemma}\label{lem:signatures}
	Let $G$ be a constant spinal group and let $d \in {\Comp(D)}$ be a compatible directed element. Let $rR' \in \spl(d)$ be a forking point with signature $(k, k')$. Then $d^k \in D$ is compatible and $r^k R'$ is a forking point of $d^k$ with signature $(1, l)$ for some $l \in (1, m(r))$.
\end{lemma}

\begin{proof}
	Let $r \in R$ and $s \in R'$. Then
	\[
		d^k|_{\underline{rs}} = (d|_{\underline{rs}})^k \equiv_{R'} (d|_{\underline{r}})^k = d^k|_{\underline{r}},
	\]
	so $d^k$ is compatible. Now let $(k, {l_2, \dots, l_{\ell}})$ be a directed path starting and ending at~$rR'$ in $\mathcal{R}^{\mathrm{ab}}(d)$ for some $k$ that is (multiplicatively) invertible modulo $\ord(d)$. Denote the vertices along {the} path by $(rR', s_2R', \dots, s_\ell R')$, thus $d|_{\underline{r^k}} \equiv_{R'} s_2$, $d|_{\underline{s_i^{l_i}}} \equiv_{R'} s_{i+1}$ for $i \in \{2, \dots, \ell-1\}$ and $d|_{\underline{s_\ell^{l_\ell}}} \equiv_{R'} r$.
	Let $k^{-1}$ be the inverse of $k$ in $\Z/\ord(d)\Z$. Then the directed path $(1, l_2k^{-1}, \dots, l_\ell k^{-1})$ starting at $r^kR'$ fulfils
	\begin{align*}
		d^k|_{\underline{r^k}} &= (d|_{\underline{r^k}})^k \equiv_{R'} s_2^k, \\
		d^k|_{\underline{(s_i^k)^{e_ik^{-1}}}} &= (d|_{\underline{s_i^{e_i}}})^k \equiv_{R'} (s_{i+1})^k \quad \text{ for } i \in [2, \ell){,}\text{ and } \\
		d^k|_{\underline{(s_\ell^k)^{e_\ell k^{-1}}}} &= (d|_{\underline{s_\ell^{e_\ell}}})^k \equiv_{R'} r^k.
	\end{align*}
	Analogously, the second directed path witnessing that $rR'$ is a forking point of $d$ yields a {second} directed path from $r^kR'$ to itself starting with an edge labelled $k'k^{-1} \neq 1$ of length equal to the {length} of the original {second} path.
\end{proof}

\begin{definition}
	Let $G$ be a group and $\mathcal{G} = (G, E)$ a directed graph with vertex set equal to the underlying set of the group $G$. Let $V \subseteq G$. The \emph{reach-closure} $\overline{\reach}(V)$ is the minimal subgroup of $G$ containing $V$ which is closed under the reach {operation within the graph} $\mathcal{G}$, i.e.\ { such that} $\reach(\overline{\reach}(V)) \subseteq \overline{\reach}(V)$.
\end{definition}

{In the remainder of the paper, we will exclusively consider directed graphs of type $\mathcal{R}^{\mathrm{ab}}(d)$ for some directed element $d$. If clarification is necessary to understand with respect to which graph the reach closure is taken, we write $\overline{\reach}_d$ to indicate that it is computed with respect to the graph associated to $d$.}

If the group $G$ is finite, the reach closure of any set $V$ can be computed by alternatingly passing to the subgroup generated by and the reach of the given subset, until one alternates between the same subset and group, which necessarily are equal to $\reach(\overline{\reach}(V))$ and $\overline{\reach}(V)$.

We shall consider the reach closure{s} of forking points, i.e. $\overline{\reach}(rR')$ for $rR' \in \spl(d)$. This is a subgroup of $R/R'$ and not of $R$; since we are rather interested in the set of possible lifts of $\overline{\reach}(rR')$, we make a final definition.

\begin{definition}\label{def:filling}
	Let $G = \langle R \cup D \rangle$ be a constant spinal group. Let $d \in {\Comp(D)}$ be a compatible directed element. We say that $d$ is \emph{filling} if there exists some $rR' \in \spl(d)$ such that
	\[
		\bigcap_{\lambda \colon \overline{\reach}(rR') \to R'} \langle t|_{\underline{s\lambda(s)}} \mid sR' \in \overline{\reach}_{d}(rR') \rangle = R
	\]
	for all $t \in d^D$, where the intersection above is taken over all maps $\lambda$ from $\overline{\reach}(rR')$ to $R'$. We say a subset $T \subseteq D$ is filling if every element $t \in T$ is filling.
\end{definition}

Note that if an element $d$ is filling, also its conjugacy class in $D$ if filling. This is a simple consequence of $\spl(d) = \spl(d^t)$ for all $t \in D$, which is immediate from the fact that $\mathcal{R}^{\mathrm{ab}}(d) =\mathcal{R}^{\mathrm{ab}}(d^t)$.

\begin{figure}
	\centering
	\begin{tikzpicture}[font=\footnotesize, scale=0.8]
		\node (id)	[circle, fill=black, inner sep=0.05cm, label=below:$\id$] at (6,2) {};
		\node (a1)	[circle, fill=black, inner sep=0.05cm, label=below:$  a$] at (2,2) {};
		\node (a2)	[circle, fill=black, inner sep=0.05cm, label=right:$a^2$] at (4,1) {};
		\node (a3)	[circle, fill=black, inner sep=0.05cm, label=right:$a^3$] at (4,3) {};

		\path[->] (a1) edge [loop left]					node			{1}	()
				  (a1) edge [bend left]					node 	[above] {2}	(a2)
				  	   edge [bend right]				node	[below]	{3}	(a2)
			      (a2) edge [loop below]				node			{1}	()
				  (a3) edge [bend left]					node 	[right] {1}	(a2)
				  	   edge [bend right]				node	[right]	{2}	(a2)
					   edge [bend right]				node	[above]	{3}	(a1);
		
		\path (7,0) edge node {} (7,3.5);
		
		\node (id2)	[circle, fill=black, inner sep=0.05cm, label=below:$\id$] at (8,2) {};
		\node (a12)	[circle, fill=black, inner sep=0.10cm, label=below:$  a$] at (10,3) {};
		\node (a22)	[circle, fill=black, inner sep=0.05cm, label=right:$a^2$] at (12,2) {};
		\node (a32)	[circle, fill=black, inner sep=0.05cm, label=right:$a^3$] at (10,1) {};

		\path[->] (a12) edge [loop above]					node			{1, 2}	()
				  (a12) edge 								node 	[above] {3}	(id2)
			      (a22) edge 								node	[above]		{1}	(a12)
				  (a32) edge 								node 	[below] {1}	(id2)
				  	    edge [bend left]					node	[right]	{2}	(a12)
				  	    edge [bend right]					node	[right]	{3}	(a12);
	\end{tikzpicture}
	\caption{
		The graph $\mathcal{R}^{\mathrm{ab}}(b)$ associated to the (primary) GGS group $\langle \{a, b\} \rangle$ with defining vector $(1, 2, 2) \in (\Z/4\Z)^3$ (left) and $(1, 1, 0) \in (\Z/4\Z)^3$ (right). The only forking point is the vertex $a$ in the second graph, whose reach is $\{\id, a\}$ and whose reach closure is the full rooted group.
	}
	\label{fig:eg pseudoggs}
\end{figure}
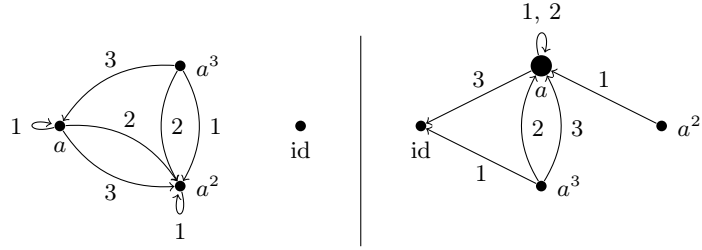


\subsection{Delta and theta maps} 
\label{sub:delta_and_theta_maps}

We now model the desired passage to `simpler' elements by certain maps. To do so, we define the following {maps}.

\begin{definition}
	Let $G = \langle R \cup D \rangle$ be a constant spinal group and let $d \in D$ be a directed element. The map $\Delta_d \colon G \to R$ is defined by
	\[
		\Delta_d(g) = \begin{cases}
			d|_{g({\underline{e}})} &\text{ if }g \notin \St_G(1),\\
			\id	&\text{ otherwise,}
		\end{cases}
	\]
	and the map $\theta_d \colon G' \times \{ (r, l) \mid r \in R, l \in [1, m(r)) \} \times \Z$ is defined by
	\[
		\theta_d(g, r, l, k) = [\Delta_d((rg)^{l-m(r)k}), (rg)^{l-m(r)k}|_{\underline{e}}].
	\]
\end{definition}

We shall study the behaviour of the passing from a pair
\[
	(r, g) \quad \text{ to the pair } \quad (\Delta_d((rg)^{l-m(r)k})), \theta_d(g,r,l,k))
\]
for fittingly chosen $l$ and $k$. The roles of $l$ and $k$ seem to be somewhat doubled and warranting an explanation. The number $l$ corresponds to (the label of an) edge in~$\mathcal{R}^{\mathrm{ab}}(d)$, thus to the passage from rooted to rooted element, while the influence of $k$ is undetectable in $R/R'$; it modifies only the non-rooted part. {With the directed graph in mind, it is thus useful to single out $l$.}

Using the definition above, we immediately arrive at the following immediate but central observation, that we record as a lemma.

\begin{lemma}\label{lem:theta in M0}
	Let $G = \langle R \cup D \rangle$ be a constant spinal group and let $M \leq G$ be a subgroup containing a compatible directed element $d \in D$. Let $k \in \Z$, let $r \in R$, let $l \in [1,m(r))$ and $g \in G'$ such that $rg \in M$. Then 
	\[
		\Delta_d((rg)^{l-m(r)k}) \theta_d(g, r, l, k) \in M_e
	\]
	Furthermore,
	\[
		\Delta_d((rg)^{l-m(r)k}) \equiv_{R'} \Delta_d(r^l),
	\]
	i.e. $\Delta_d((rg)^{l-m(r)k})R'$ is the unique vertex in $\mathcal{R}^{\mathrm{ab}}(d)$ that is reached from $rR'$ through an edge labeled $l$.
\end{lemma}

\begin{proof}
	Write $t = {l-m(r)k}$. First note that if $(rg)^t \in \St(1)$, both $\Delta_d((rg)^t)$ and $\theta_d(g, r, l, k)$ are trivial, hence the first statement holds. Assume otherwise. Then the first claim follows from a direct calculation, using that $d^{(rg)^t} \in {\st}_M({\underline{e}})$
	\begin{align*}
		M_{\underline{e}} \ni d^{(rg)^t}|_{\underline{e}} &= (d|_{(rg)^t({\underline{e}})})^{(rg)^t|_{\underline{e}}}\\&= d|_{(rg)^t({\underline{e}})} [d|_{(rg)^t({\underline{e}})}, (rg)^t|_{\underline{e}}] = \Delta_d((rg)^t) \theta_d(g, r, l, k).
	\end{align*}
	To conclude, recall that $d$ is compatible, so $\Delta_d(g) \equiv_{R'} \Delta_d(g y)$ for all $y \in G'$. Therefore, if $(rg)^{l-m(r)k} \equiv_{G'} r^{l-m(r)k}\equiv_{R'} r^l$, then there exists some $s\in R'$ such that $\Delta_d((rg)^{l-m(r)k}) = \Delta_d(r^l)s$.
\end{proof}

It is difficult to {monitor these} mappings. A particular problem is that the rooted subgroup and the commutator subgroup of $G$ intersect in $R'$ due to \cref{thm:factor commutator groups of constant spinal groups}. Thus all control over the rooted part of the product $rg$ is only up to the subgroup $R'$.


\subsection{Contraction under the maps $\theta_d$} 
\label{sub:theta_maps}

We now establish conditions under which the syllable length of an element $g$ decreases under the application of a map $g \mapsto \theta_d(g, r, l, k)$ for suitable choice of $r \in R$, $l \in [1, m(r))$ and $k \in \Z$. Naturally, we {make heavy use of} \cref{lem:contraction lemma}. We establish some further terminology.

\begin{definition}
	Let $G = \langle R \cup D \rangle$ be a constant spinal group. Let $Y \subseteq X$. An element $g \in G$ is called \emph{saturated in $Y$} if $g|_x \in R$ for all $x \in X \smallsetminus Y$. It is called \emph{properly saturated in $Y$} if $g|_x \in R$ if and only if $x \notin Y$.
	
	Furthermore, an element of even syllable length is called \emph{split-saturated in $Y_0$ and $Y_1$} for some subsets $Y_0, Y_1 \subseteq X$ if $g$ is saturated in $Y_0 \cup Y_1$ and $\sum_{y \in Y_i} \syl(g|_y) = \syl(g)/2$ for $i \in \{0, 1\}$.
\end{definition}

Evidently every element is saturated in $X$. If an element is given in minimal syllable form $g = r_0 \prod_{i = 1}^n d_i^{r_i}$, then it is properly saturated in $\{ r_i^{-1} \mid i \in [1,n] \}$. If $g \in G$ is properly saturated in $Y$ and $r \in R$, then $g^r$ is properly saturated in $r^{-1}Y$.

\begin{lemma}\label{lem: split saturated}
	Let $G = \langle R \cup D \rangle$ be a constant spinal group and let $g \in G$ be given in minimal syllable form $g = s \prod_{i = 1}^n d_i^{r_i}$. If $g$ is split-saturated in {a singleton} $\{r\} {\subset R}$ and some set $Y \subseteq X$, then either
	\begin{align*}
		r_i^{-1} &= r \quad\text{ if and only if {$i$ is even}, or}\\
		r_i^{-1} &= r \quad\text{ if and only if {$i$ is odd.}}
	\end{align*}
\end{lemma}

\begin{proof}
	Since $g$ is split-saturated in $\{r\}$ and $Y$, half of the $i \in \{1, \dots, n\}$  must satisfy $d_i^{r_i}|_{\underline{r}}\in D$, hence $r_i = r^{-1}$. On the other hand, $r_i \neq r_{i+1}$ since any minimal syllable form is reduced. Therefore precisely every second syllable satisfies $d_i^{r_i}|_{\underline{r}}\in D$ and in particular $r_i = r^{-1}$.
\end{proof}

\begin{lemma}\label{lem:basic theta red}
	Let $G = \langle R \cup D \rangle$ be a constant spinal group, let $d \in D$, let $r \in R$ and let $g \in G'$. Write $sh$ for the decomposition of $g$ into a product of element $s \in R'$ and $h \in \St_{G}(1)$. Denote by $k$ the quotient $\ord_R(rs)/m(r)${, such that $(rs)^{km(r)} = e$}. Let $l \in [1, \dots, m(r))$. Then
	\[
		\min\{ \syl(\theta_d(g, r, l, 0)), \syl(\theta_d(g, r, l, k))\} \leq \syl(g).
	\]
	Furthermore, if the above is an equality, then
	\begin{enumerate}
		\item the element $g$ is split-saturated in the two disjoint parts
		\[
			\{ 1, rs, (rs)^2, \dots, (rs)^{l-1}\} \quad\text{and}\quad \{ (rs)^{l}, (rs)^{l+1}, \dots, (rs)^{-1}\}
		\]
		of the subgroup $\langle rs \rangle$, and
		\item $\syl(\theta_d(g, r, l, v)) = 2 \syl((rg)^{l-m(r)v}|_{\underline{e}})$ for {all} $v \in \{0, k\}$.
	\end{enumerate}
	In particular, in this case the element $g$ {necessarily has} even syllable length.
\end{lemma}

\begin{proof}
	We show that
	\[
		\syl(\theta_d(g, r, l, 0)) + \syl(\theta_d(g, r, l, k)) \leq 2\syl(g),
	\]
	from which one readily derives the inequality in the lemma.
	
	Let $v \in \{0, k\}$ and write $r' = d|_{(rg)^{l-m(r)v}({\underline{e}})} = d|_{\underline{(rs)^{l-m(r)v}}} \in R$. Unraveling the definition, we find
	\begin{align*}
		\theta_d(g, r, l, v) &= [\Delta_d((rg)^{l-m(r)v}),(rg)^{l-m(r)v}|_{\underline{e}}] = [d|_{(rg)^{l-m(r)v}({\underline{e}})},(rg)^{l-m(r)v}|_{\underline{e}}]\\ &= [r', (rg)^{l-m(r)v}|_{\underline{e}}] = (((rg)^{l-m(r)v}|_{\underline{e}})^{-1})^{r'} (rg)^{l-m(r)v}|_{\underline{e}},
	\end{align*}
	providing the inequality 
	\begin{align*}
		\syl(\theta_d(g, r, l, v)) \leq 2 \syl((rg)^{l-m(r)v}|_{\underline{e}}).
	\end{align*}
	Furthermore,
	\begin{align*}
		(rg)^l|_{\underline{e}} &= (rg)|_{(rg)^{l-1}({\underline{e}})}\cdots(rg)|_{(rg)({\underline{e}})}(rg)|_{\underline{e}}
		= g|_{\underline{(rs)^{l-1}}} \dots g|_{\underline{rs}} g|_{\underline{e}},
	\end{align*}
	and, using the above and the fact that $m(r)k = \ord_R(rs)$,
	\begin{align*}
		(rg)^{l-m(r)k}|_{\underline{e}} &= ((rg)^{-1})^{m(r)k-l}|_{\underline{e}} = ((rs)^{-1}(h^{-1})^{(rs)^{-1}})^{m(r)k-l}|_{\underline{e}}\\
		&= (h^{-1})^{(rs)^{-1}}|_{\underline{(rs)^{l-m(r)k+1}}} \dots (h^{-1})^{(rs)^{-1}}|_{\underline{(rs)^{-1}}}(h^{-1})^{(rs)^{-1}}|_{\underline{e}}\\
		&= (h|_{\underline{(rs)^{l}}})^{-1} \dots (h|_{\underline{(rs)^{-2}}})^{-1}(h|_{\underline{(rs)^{-1}}})^{-1}\\
		&= (g|_{\underline{(rs)^{-1}}}g|_{\underline{(rs)^{-2}}} \dots g|_{\underline{(rs)^{l+1}}}g|_{\underline{(rs)^{l}}})^{-1}.
	\end{align*}
	Now we find
	\begin{align*}
		\syl(\theta_d(g, r, l, 0)) + \syl(\theta_d(g, r, l, k)) &\leq 2 (\syl((rg)^{l}|_{\underline{e}}) + \syl((rg)^{l-m(r)k}|_{\underline{e}}))\\
		&\leq 2\left(\sum_{i = 0}^{l-1} \syl(g|_{\underline{(rs)^i}}) + \sum_{i = l}^{m(r)k-1} \syl(g|_{\underline{(rs)^i}})\right)\\
		&\leq 2 \syl(g),
	\end{align*}
	where we have used \cref{lem:contraction lemma}(3) in the last step.
	
	If we have
	\[
		\min\{ \syl(\theta_d(g, r, l, 0)), \syl(\theta_d(g, r, l, k))\} = \syl(g),
	\]
	we see that all inequalities in the previous computation are in fact equalities. Thus both $\sum_{i = 0}^{l-1} \syl(g|_{(rs)^i})$ and $\sum_{i = l}^{m(r)k-1} \syl(g|_{(rs)^{i}})$ are equal to $\syl(g)/2$, which implies that $g$ is split-saturated in the two parts of the cyclic subgroup $\langle rs \rangle$ described in the statement of the lemma. Looking at the first inequality, we find that
	\[
		\syl(\theta_d(g, r, l, v)) = 2 \syl((rg)^{l-m(r)v}|_{\underline{e}})
	\]
	for $v \in \{0, k\}$.
\end{proof}

\begin{lemma}\label{lem:theta red}
	Let $G = \langle R \cup D \rangle$ be a constant spinal group. Let $h \in G$ be an element of syllable length at least $2$, let $r \in R$ be a non-trivial element, let $d \in D$ and write $g = [r, h]$. Assume that $\syl(g) = 2 \syl(h)$. Let $m(r)k$ be the order of $rh$ in $G/\St_G(1)$ and let $l \in (1, m(r))$. Then we find
	\[
		\min\{ \syl(\theta_d(g, r, u, v)) \mid u \in \{1, l\}, v \in \{0,k\}\} < \syl(g).
	\]
\end{lemma}

\begin{proof}
	First note that the existence of $l \in (1, m(r))$ implies $m(r) > 2$, hence the order of {$rR'$} in $R/R'$ is at least $3$.
	
	Let $h = s \prod_{i = 1}^n d_i^{r_i}$ be a minimal syllable form for $h$, hence $\syl(g) = 2n$. Write $\tilde{h} = s^{-1}h \in \St_G(1)$. Consider that since $r$ is non-trivial,
	\[
		g = [r, s\tilde{h}] = [r, s] (\tilde{h}^{-1})^{r^s} \tilde{h} = [r, s] \prod_{i = 1}^n (d_{n-i}^{-1})^{r_ir^s} \prod_{i = 1}^n d_i^{r_i}
	\]
	is a syllable form of length $2n$, hence a minimal syllable form.
	
	If one of the values $\syl(\theta_d(g, r, u, v))$ exceeds $\syl(g)$, the statement of the lemma follows immediately by \cref{lem:basic theta red}. Thus, assume for contradiction that all four values of $\syl(\theta_d(g, r, u, v))$ are equal to $\syl(g)$.
	
	Then, using the second part of \cref{lem:basic theta red}, the element $g$ is both split-saturated in $\{ {\underline{e}} \}$ and {$\{\underline{r[r,s]} = \underline{r^s}, \underline{(r^s)^2}, \dots, \underline{(r^s)^{-1}}\}$}, using $u = 1$, as well as split-saturated in {$\{\underline{e}, \dots, \underline{(r^s)^{l-1}}\}$ and $\{\underline{(r^s)^l}, \dots, \underline{(r^s)^{-1}}\}$}, using $u = l$.
	
	In view of \cref{lem: split saturated}, we have that precisely every second syllable in the syllable form for $g$ above is contained in $D$. Since the given minimal syllable form of $g$ ends in the given (minimal) syllable form of $\tilde{h}$, the same holds for $\tilde{h}$, i.e.\ $\tilde{h}$ starts in one of the following ways:
	\[
		\tilde{h} = \begin{cases}
			d_1^{r_1} d_2 \dot{h} &\text{ with } d_1, d_2, r_1 \neq 1, \dot{h} \in \St_G(1) \text{ or}\\
			d_1 d_2^{r_2} \dot{h} &\text{ with } d_1, d_2, r_2 \neq 1, \dot{h} \in \St_G(1).
		\end{cases}
	\]
	Here we use that $\syl(\tilde{h}) = \syl(h) \geq 2$. Now we compute
	\[
		g = \begin{cases}
			[r, s] (\dot{h}^{-1})^{r^s} (d_2^{-1})^{r^s} (d_1^{-1})^{r_1r^s} d_1^{r_1} d_2 \dot{h} &\text{ in the first case, or }\\
			[r, s] (\dot{h}^{-1})^{r^s} (d_2^{-1})^{r_2r^s} (d_1^{-1})^{r^s} d_1 d_2^{r_2} \dot{h}
			&\text{ in the second case.}
		\end{cases}
	\]
	Since every second syllable of $g$ must be in $D$ and $r$ is not trivial, this implies $r_1 = (r^s)^{-1}$ or $r_2 = (r^s)^{-1}$, respectively. In both cases, $g$ is properly saturated in a set $I$ that contains {$\{\underline{e}, \underline{r^s}, \underline{(r^s)^{-1}}\}$}. Since $r^s$ has order greater than $2$, the set $I$ has at least $3$ members. But being both split-saturated in $\{{\underline{e}}\}$ and {$\{\underline{r^s}, \underline{(r^s)^2}, \dots, \underline{(r^s)^{-1}}\}$} and being properly saturated in $I$ implies, for all $l \in (1, m(r))$,
	\[
		\sum_{i = 0}^{l-1} \syl(g|_{\underline{(r^s)^i}}) \geq \syl(g|_{\underline{e}}) + \syl(g|_{\underline{r^s}}) > \syl(g|_{\underline{e}}) = \syl(g)/2.
	\]
	Thus it is impossible for $g$ to be split-saturated in the sets {$\{\underline{e}, \underline{r^s}, \dots, \underline{(r^s)^{l-1}}\}$ and $\{\underline{(r^s)^l}, \dots \underline{(r^s)^{-1}}\}$}. So $\syl(\theta_d(g, r, l, k)) < \syl(g)$.
\end{proof}


\subsection{Reduction along paths} 
\label{sub:reduction_along_paths}

The passage modelled by $\Delta_d$ and $\theta_d$ corresponds to the change of vertex $rR'$ to another vertex $\Delta((rg)^{l-m(r)k})$ that is connected to $rR'$ by a directed edge labeled $l$. We now extend the definitions of $\Delta_d$ and $\theta_d$ by considering a directed path $\mathfrak{p}$ of length $\ell$ within $\mathcal{R}^{\mathrm{ab}}(d)$ instead of the length $1$ path described by the label $l \in [1, m(r))$. Let $t \in R/R'$ be the initial vertex of $\mathfrak{p}$, let $r \in tR'$ and let $g \in G'$. Furthermore, let $\kappa = (k_i)_{i \in [1, \ell]} \in \Z^{[1, \ell]}$. Define
\[
	\Delta_d(g, r, \mathfrak{p}, \kappa) = r \quad\text{and}\quad \theta_d(g, r, \mathfrak{p}, \kappa) = g
\]
in case of the trivial path of length $0$. If $\mathfrak{p}$ is any other path, denote by $\mathfrak{q}$ the sub-path consisting of all but the last edge of $\mathfrak{p}$, denote by $l$ the label of the last edge of $\mathfrak{p}$, denote by $m = m(s)$ the order corresponding of the final vertex $sR'$ of $\mathfrak{q}$, and denote by $\delta = (k_i)_{i \in [1, \ell)}$. Then set
\begin{align*}
	\Delta_d(g, r, \mathfrak{p}, \kappa) &= \Delta_d((\Delta_d(g, r, \mathfrak{q}, \delta) \theta_d(g, r, \mathfrak{q}, \delta))^{l - m k_{\ell}})\\
	\theta_d(g, r, \mathfrak{p}, \kappa) &= \theta_d( \theta_d(g, r, \mathfrak{q}, \delta), \Delta_d(g, r, \mathfrak{q}, \delta),l, k_{\ell})\\
	&=[\Delta_d(g, r, \mathfrak{p}, \kappa), (\Delta_d(g, r, \mathfrak{q}, \delta)\theta_d(g, r, \mathfrak{q}, \delta))^{l - m k_{\ell}}|_{\underline{e}}]
\end{align*}
Evidently, if $\mathfrak{p}$ is a path of length $1$, this corresponds to the original maps $\theta_d$ and $\Delta_d$ for suitable values of $r$, $l$ and $k$.

\begin{lemma}\label{lem:path lemma}
	Let $G = \langle R \cup D \rangle$ be a constant spinal group. Let $d \in {\Comp(D)}$ be compatible, let $\mathfrak{p}$ be a directed path of length $\ell$ within $\mathcal{R}^{\mathrm{ab}}(d)$ from $t_1 \in R/R'$ to $t_2 \in R/R'$. Let $r \in t_1R'$ and let $g \in G'$. Let $M \leq G$ be a subgroup containing the product $rg$. Then there exists $\kappa = (k_i)_{i \in [1, \ell]} \in \Z^{[1, \ell]}$ such that
	\[
		\Delta_d(g, r, \mathfrak{p}, \kappa) \theta_d(g, r, \mathfrak{p}, \kappa) \in M_{{\underline{e}}^\ell}
	\]
	and $\syl(\theta_d(g, r, \mathfrak{p}, \kappa)) \leq \syl(g)$, and $\Delta_d(g, r, \mathfrak{p}, \kappa) \equiv_{R'} t_2$.
	
	Furthermore, if $\syl(\theta_d(g, r, \mathfrak{p}, \kappa)) = \syl(g)$, we may assume $\theta_d(g, r, \mathfrak{p}, \kappa)=[\tilde{r}, h]$ and $\syl(\theta_d(g, r, \mathfrak{p}, \kappa)) = 2\syl(h)$ for some $\tilde{r}\in R$ and $h\in G$.
\end{lemma}

\begin{proof}
	The first statement follows by applying \cref{lem:theta in M0} and \cref{lem:basic theta red} repeatedly, the choice of  is $\kappa$ dictated by the second lemma. Consult the second part of \cref{lem:theta in M0}, to see that $\Delta_d(g, r, \mathfrak{p}, \epsilon) \in t_2R'$. The last part of the lemma follows by the definition of $ \theta_d(g, r, \mathfrak{p}, \kappa)$ and \cref{lem:basic theta red}(2).
\end{proof}

Furthermore, we are able to concisely {reformulate} the {results} of the previous subsection in a single lemma.

\begin{lemma}\label{lem:theta red comb}
	Let $g \in G$, let $r \in R$ be non-trivial, let $d \in {\Comp(D)}$ {be compatible}, let $u \in {[m(r))}$ and let $l \in (1, m(\Delta_d((rg)^{u-m(r)k})))$. Assume that $\syl(g) > 2$. Then there exist $k, k_1, k_2 \in \Z$ such that either
	\begin{align*}
		\syl(\theta_d(g, r, (u, 1), (k, k_1)) &< \syl(g) \quad\text{or}\\
		\syl(\theta_d(g, r, (u, l), (k, k_2)) &< \syl(g).
	\end{align*}
\end{lemma}
\begin{proof}
	In case $\syl(g) = 3$, by \cref{lem:basic theta red} there exists $k \in \Z$ such that
	\[
		\syl(\theta_d(g, r, u, k)) < \syl(g)
	\]
	for $i \in \{1, 2\}$. Using \cref{lem:basic theta red} again, we find $k_i$ for $i \in \{1, 2\}$ such that both inequalities hold.
	Now assume $\syl(g) > 3$. Using \cref{lem:basic theta red} we find $k\in \Z$ such that
	\[
		\syl(\theta_d(g, r, u, k)) \leq \syl(g).
	\]
	If the inequality is strict, we continue as before; otherwise, we find that
	\[
		\syl(\theta_d(g, r, u, k)) = \syl([\Delta_d((rg)^{u-m(r)k}), (rg)^{u-m(r)k}]|_{\underline{e}}) = \syl(g) > 3
	\]
	implies that $\syl((rg)^{u-m(r)k}]|_{\underline{e}}) \geq 2$. Also, by \cref{lem:basic theta red}, the inequality being non-strict implies $\syl(\theta_d(g, r, u, k)) = 2 \syl((rg)^{u-m(r)k}]|_{\underline{e}})$. Thus we may apply \cref{lem:theta red} and find $k_1, k_2 \in \Z$ such that the statement holds.
\end{proof}

The iterative definition of the extended definitions of $\Delta_d$ and $\theta_d$ makes it easy to observe that given two paths $\mathfrak{p}$ and $\mathfrak{q}$ of length $\ell'$ and $\ell''$, respectively, such that the composition $\mathfrak{p}\circ\mathfrak{q}$ is well-defined --- i.e.\ such that the last vertex of $\mathfrak{p}$ coincides with the first of $\mathfrak{q}$ --- and given two tuples $\kappa' \in \Z^{[1, \ell']}$ and $\kappa'' \in \Z^{[1, \ell'']}$ and their concatenation $\kappa' \circ \kappa'' = (k'_1, \dots, k'_{\ell'}, k''_1, \dots, k''_{\ell''}) \in \Z^{[1, \ell'+\ell'']}$, we find
\[
	\Delta_d(g, r, \mathfrak{p}\circ\mathfrak{q}, \kappa' \circ \kappa'') = \Delta_d(\theta_d(r,g, \mathfrak{p}, \kappa'), \Delta(r,g, \mathfrak{p}, \kappa'), \mathfrak{q}, \kappa'')
\]
and
\[
	\theta_d(g, r, \mathfrak{p}\circ\mathfrak{q}, \kappa' \circ \kappa'') = \theta_d(\theta_d(g, r, \mathfrak{p}, \kappa'), \Delta_d(g, r, \mathfrak{p}, \kappa'), \mathfrak{q}, \kappa'').
\]

\begin{lemma}\label{lem:red path}
	Let $g \in G$ of syllable length at least $3$, let {$d \in \Comp(D)$ be a compatible element}, and let $r \in R$ such that there exist
	\begin{enumerate}
		\item a non-trivial directed path $\mathfrak{r}$ from $rR'$ to itself,
		\item a directed path $\mathfrak{p}$ starting at $rR'$ with first label $1$,
		\item a directed path $\mathfrak{q}$ starting at $rR'$ with first label $l \neq 1$.
	\end{enumerate}
	Then there exist $\kappa_1 \in \Z^{|\mathfrak{r}| + |\mathfrak{p}|}$ and $\kappa_2 \in \Z^{|\mathfrak{r}| + |\mathfrak{q}|}$ such that either
	\begin{align*}
		\syl(\theta_d(g, r, \mathfrak{r} \circ \mathfrak{p}, \kappa_1)) &< \syl(g) \quad\text{or}\\
		\syl(\theta_d(g, r, \mathfrak{r} \circ \mathfrak{q}, \kappa_2)) &< \syl(g).
	\end{align*}
\end{lemma}

\begin{proof}
	Write $\mathfrak{r} = \mathfrak{a} \circ (u)$, $\mathfrak{p} = (1) \circ \mathfrak{o}$ and $\mathfrak{q} = (l) \circ \mathfrak{e}$ for some potentially trivial directed paths $\mathfrak{a}, \mathfrak{o}$ and $\mathfrak{e}$. Hence
	\begin{align*}
		\mathfrak{r} \circ \mathfrak{p} &= \mathfrak{a} \circ (u, 1) \circ \mathfrak{o} \quad\text{and}\\
		\mathfrak{r} \circ \mathfrak{q} &= \mathfrak{a} \circ (u, l) \circ \mathfrak{e}.
	\end{align*}
	By \cref{lem:path lemma}, we find $\epsilon \in \Z^{\ell-1}$ such that
	\[
		\syl(\theta_d(g, r, \mathfrak{a}, \epsilon)) \leq \syl(g).
	\]
	Now by \cref{lem:theta red comb}, we find $(k, k_1)$ and $(k, k_2)$ such that either
	\begin{align*}
		\syl(\theta_d(g, r, \mathfrak{a} &\circ(u, 1), \epsilon \circ (k, k_1))\\
		&= \syl(\theta_d(\theta_d(g, r, \mathfrak{a}, \epsilon), \Delta_d(g, r, \mathfrak{a}, \epsilon), (u, 1), (k, k_1))) < \syl(g) \quad\text{ or }\\
		\syl(\theta_d(g, r, \mathfrak{a} &\circ(u, l), \epsilon \circ (k, k_2))\\
		&= \syl(\theta_d(\theta_d(g, r, \mathfrak{a}, \epsilon), \Delta_d(g, r, \mathfrak{a}, \epsilon), (u, l), (k, k_2))) < \syl(g).
	\end{align*}
	Using \cref{lem:path lemma}, we find $\delta \in \Z^{\ell-1}$ and $\eta \in \Z^{\ell'-1}$ such that
	\begin{align*}
		\syl(\theta_d(&g, r, \mathfrak{a} \circ(u, 1) \circ \mathfrak{o}, \epsilon \circ (k, k_1) \circ \delta)\\
		&= \syl(\theta_d(\theta_d(g, r, \mathfrak{a} \circ(u, 1), \epsilon \circ (k, k_1)), \Delta_d(g, r, \mathfrak{a} \circ (u, 1), \epsilon \circ (k, k_1)), \mathfrak{o}, \delta))\\
		&\leq \syl(\theta_d(g, r, \mathfrak{a} \circ(u, 1), \epsilon \circ (k, k_1))
	\end{align*}
	and
	\begin{align*}
		\syl(\theta_d(&g, r, \mathfrak{a} \circ(u, l) \circ \mathfrak{e}, \epsilon \circ (k, k_2) \circ \eta)\\
		&= \syl(\theta_d(\theta_d(g, r, \mathfrak{a} \circ(u, l), \epsilon \circ (k, k_2)), \Delta_d(g, r, \mathfrak{a} \circ (u, l), \epsilon \circ (k, k_2)), \mathfrak{e}, \eta))\\
		&\leq \syl(\theta_d(g, r, \mathfrak{a} \circ(u, l), \epsilon \circ (k, k_2)).
	\end{align*}
	Thus we may choose $\kappa_1 = \epsilon \circ (k, k_1) \circ \delta$ and $\kappa_2 = \epsilon \circ (k, k_2) \circ \eta$ to obtain the statement.
\end{proof}

\begin{lemma}\label{lem:core from once}
	Let $G = \langle R \cup D \rangle$ be a constant spinal group. Let $M \leq G$ be a subgroup containing both a compatible element $d \in {\Comp(D)}$ and {a} rooted element $r \in R$ such that $rR' \in \spl(d)$ {is a forking point}. Then $rR' \in \pi^{\mathrm{ab}}_{R}(\core(M) \cap R)$.
\end{lemma}

\begin{proof}
	Since $rR'$ {is a forking point,} there exist two directed paths of coprime lengths $\ell_{\mathfrak{p}}, \ell_{\mathfrak{q}}$ from $rR'$ to itself. Consequently, using the fact that every integer greater than $n = \ell_{\mathfrak{p}}\ell_{\mathfrak{q}} - \ell_{\mathfrak{p}} - \ell_{\mathfrak{q}}$ can be represented as a positive integer combination of $\ell_{\mathfrak{p}}$ and $\ell_{\mathfrak{q}}$, we find a path $\mathfrak{r}_N$ from $rR'$ to itself of length $N$ for all $N > n$. Now by \cref{lem:path lemma}, there exists some $\kappa_N \in \Z^{N}$ such that
	\begin{align*}
		\Delta_d(\id, r, \mathfrak{r}_N, \kappa_N) \theta_d(\id, r, \mathfrak{r}_N, \kappa_N) &\in M_{{\underline{e}}^N},\\
		\syl(\theta_d(\id, r, \mathfrak{r}_N, \kappa_N)) &\leq \syl(\id) = 0, \quad\text{ and}\\
		\Delta_d(\id, r, \mathfrak{r}_N, \kappa_N)&\equiv_{R'}r.
	\end{align*}
	Using the fact that $\theta_d(\id, r, \mathfrak{r}_N, \kappa_N) \in G'$, we find
	\[
		rR' \in \pi^{\mathrm{ab}}_{R}(M_{{\underline{e}}^N} \cap R).
	\]	
	Since $N$ {is an arbitrary integer greater than $n$}, the same is also true for $\core(M)$ instead of~$M_{{\underline{e}}^N}$.
\end{proof}

\begin{lemma}\label{lem:forks in core}
	Let $G = \langle R \cup D \rangle$ be a constant spinal group. Let $d \in {\Comp(D)}$ be a compatible element and $rR' \in \spl(d)$. Let $M \leq G$ be a prodense subgroup of $G$ containing $d$. Then there exists an element $t \in D$ such that
	\[
		rR' \in \pi^{\mathrm{ab}}_{R}(\core(M^{t}) \cap R).
	\]
\end{lemma}

\begin{proof}
	Since $M$ is prodense and $G'$ is a normal subgroup of finite index in $G$, the intersection $rG' \cap M$ is non-empty, hence we find $g \in G'$ such that $rg \in M$.
	
	The element $rR'$ is contained in $\spl(d)$. In view of \cref{lem:signatures}, we may assume that $rR'$ has signature $(1, l)$ for some $l \in (1, m(r))$, i.e.\ there are at least two distinct directed paths $\mathfrak{p}$ and $\mathfrak{q}$ of coprime lengths $|\mathfrak{p}| = \ell$ and $|\mathfrak{q}| = \ell'$ leading from $rR'$ to itself; path $\mathfrak{p}$ starts with an edge labeled by $1$ and the second path $\mathfrak{q}$ starts with an edge labeled by $l$.
	
	Assume that $\syl(g) \geq 3$. By \cref{lem:red path}, we find that either $\theta_d(g, r, \mathfrak{p} \circ \mathfrak{p}, \kappa)$ or $\theta_d(g, r, \mathfrak{p} \circ \mathfrak{q}, \kappa)$ for some $\kappa$ of appropriate length has {a} shorter syllable length {than $g$}, while $\Delta_d(g, r, \mathfrak{p} \circ \mathfrak{p}, \kappa)$ {and} $\Delta_d(g, r, \mathfrak{p} \circ \mathfrak{q}, \kappa)$ {are} equivalent to $r$ modulo $R'$ and the {respective} product{s 
	\[
		\theta_d(g,r,\mathfrak{p}\circ\mathfrak{p}, \kappa)\Delta_d(g,r,\mathfrak{p}\circ\mathfrak{p}, \kappa) \quad\text{ and }\quad \theta_d(g,r,\mathfrak{p}\circ\mathfrak{q}, \kappa)\Delta_d(g,r,\mathfrak{p}\circ\mathfrak{q}, \kappa)
	\]
	are} contained in the projection of $M_{{\underline{e}}^{2\ell}}$ or $M_{{\underline{e}}^{\ell+\ell'}}$, respectively. Since the core of $M$ is equal to the core of any projection $M_{e^k}$, we may replace $g$ by its image under $\theta_d$ of shorter syllable length and $r$ by the corresponding image under $\Delta_d$. Repeating this process, we may assume that $\syl(g) < 3$.
	
	Furthermore, using the last paragraph of \cref{lem:path lemma}, we may assume that $g = [r, h]$ for $h \in G$ such that $\syl(g) = 2 \syl(h)$, or $\syl(g) \leq 1$, by replacing $g$ (yet another time) with $\theta_d(g, r, \mathfrak{p}, \kappa)$ for a suitable $\kappa$ and $r$ by the corresponding image under $\Delta_d$.

	Consider first the case that $\syl(g) \leq 1$. Then by \cref{lem:basic theta red} there exists some $k \in \Z$ such that either
	\[
		\syl(\theta_d(g, r, 1, k)) = 0 \quad \text{ or } \quad \syl(\theta_d(g, r, 1, 0)) = 0.
	\]
	Proceeding along the rest of the path $\mathfrak{p}$, we may indeed replace $g$ and $r$ with suitable images under $\theta_d$ and $\Delta_d$, respectively, allowing us to assume that $\syl(g) = 0$. Thus $g \in R \cap G' = R'$, hence $rg \in M \cap R$. By \cref{lem:core from once}, we find that $rR' \in \pi^{\mathrm{ab}}_{R}(\core(M) \cap R)$.

	Thus it remains {to consider} the case that $\syl(g) = 2$ and $g = [r, h]$ with $\syl(g) = 2 \syl(h)$. Clearly this implies $\syl(h) = 1$. Thus we write $h = s_0 t^{s_1}$ for some $s_0, s_1 \in R$ and $t \in D$. In consequence, $g = [r, s_0] (t^{-1})^{s_1r^{s_0}}t^{s_1}$.
		
	We use \cref{lem:basic theta red} to restrict our attention to the case that $s_1 \in \{e, r^{s_0}, (r^{s_0})^{-1}\}$. For a simpler notation, set $\dot{r} = r^{s_0}$. By \cref{lem:basic theta red}, either $g$ is split-saturated in $\{{\underline{e}}\}$ and {$\{ \underline{r[r, s_0]} = \underline{\dot{r}}, \dots, \underline{\dot{r}^{-1}} \}$}, or we may reduce the length of $g$ by replacing it with a shorter element $\theta_d(g, r, \mathfrak{p}, \kappa)$ for a fitting choice of $\kappa$. {The latter option} puts us in the situation of the previous case, and we are done.
	
	If $g = [r, s_0] (t^{-1})^{s_1\dot{r}}t^{s_1}$ is split-saturated in $\{{\underline{e}}\}$ and {$\{ \underline{\dot{r}}, \dots, \underline{\dot{r}^{-1}} \}$}, it necessarily follows that $s_1 = e$ or $s_1^{-1} = \dot{r}$. If $s_1 = e$, we find
	\[
		g = [r, s_0] (t^{-1})^{\dot{r}}t = r^{-1}(\dot{r})^{t},
	\]
	hence $rg = \dot{r}^{t} \in M$. Thus $M^{t^{-1}}$ contains both $d^{t^{-1}}$ and $\dot{r} \equiv_{R'} r$, whence an application of \cref{lem:core from once} yields the statement of the lemma.
	
	Now consider the case $s_1^{-1} = \dot{r}$. Using this, we find
	\[
		g = [r, s_0] t^{-1}t^{\dot{r}^{-1}},
	\]
	hence
	\[
		rg = \dot{r} t^{-1} \dot{r} t \dot{r}^{-1} = (\dot{r})^{ t \dot{r}^{-1} }.
	\]
	Consider that both $d$ and $rg$ are contained in $M$, so
	\begin{align*}
		d^{(rg)^l}|_{\underline{e}} &= d^{(\dot{r}^l)^{t\dot{r}^{- 1}}}|_{\underline{e}} = (d|_{(\dot{r}^l)^{t\dot{r}^{-1}}({\underline{e}})})^{((\dot{r}^l)^{t\dot{r}^{-1}})|_{\underline{e}}}\\
		&= (d|_{\underline{\dot{r}^l}})^{t^{-1}|_{\underline{\dot{r}^{l-1}}}t|_{\underline{\dot{r}^{-1}}}} = \Delta_{d}(\dot{r}^l)^{t^{-1}|_{\underline{\dot{r}^{l-1}}}t|_{\underline{\dot{r}^{-1}}}} \in M_{\underline{e}},
	\end{align*}
	with $\widehat{r} = \Delta_{d}(\dot{r}^l)^{t^{-1}|_{\underline{\dot{r}^{l-1}}}t|_{\underline{\dot{r}^{-1}}}} \in R$ since $l \neq 1$, and clearly $\widehat{r} \equiv_{R'} \Delta_d(\dot{r}^l)$. Following the remainder of the path $\mathfrak{q}$ from $\Delta_d(\dot{r}^l)$ to $rR'$, the statement of the lemma by a final application of \cref{lem:path lemma} and \cref{lem:core from once}.
\end{proof}

\begin{lemma}\label{lem:reach in core}
	Let $G = \langle R \cup D \rangle$ be a constant spinal group, let $d \in {\Comp(D)}$ be a compatible element, let $M$ be a prodense subgroup of $G$, and let $Y \subseteq R/R'$. If
	\[
		Y \subseteq \pi^{\mathrm{ab}}_{R}(\core(M) \cap R),
	\]
	then also
	\[
		\reach(Y) \subseteq \pi^{\mathrm{ab}}_{R}(\core(M) \cap R).
	\]
\end{lemma}

\begin{proof}
	Let $n \in \N$ be such that $Y \subseteq \pi^{\mathrm{ab}}_{R}(M_{{\underline{e}}^N} \cap R)$ for all $N > n$. Let $s \in \reach_{d}(Y)$. Then there exists a directed path $\mathfrak{p}_s$ of length $\ell$ from a vertex $y$ of $Y$ to $s$. Let $r \in M_{\underline{{e}}^N} \cap yR'$ for some $N > n$. By \cref{lem:path lemma}, we find some $\kappa \in \Z^{{[\ell]}}$ such that
	\[
		\Delta_d(\id, r, \mathfrak{p}_s, \kappa) \theta_d(\id, r, \mathfrak{p}, \kappa) \in (M_{{\underline{e}}^N})_{{\underline{e}}^\ell} = M_{{\underline{e}}^{N+\ell}}.
	\]
	Since $\syl(\theta_d(\id, r, \mathfrak{p}, \kappa)) \leq \syl(\id)$, we find $\Delta_d(\id, r, \mathfrak{p}_s, \kappa)R' = sR' \in \pi^{\mathrm{ab}}_{R}(M_{{\underline{e}}^{N+\ell}} \cap R)$ for all $N > n$. Thus $sR' \in \pi^{\mathrm{ab}}_{R}(\core(M) \cap R)$.
\end{proof}

\begin{lemma}\label{lem:special in prodense}
	Let $G = \langle R \cup D \rangle$ be a constant spinal group, let $d \in {\Comp(D)}$ be a compatible and filling element, and let $M \leq G$ be a prodense subgroup of $G$ containing $d$. Then there exists an element $t \in D$ such that
	\[
		\{ d^{t} \} \cup R \subseteq \core(M^{t}).
	\]
	In particular, a projection of $M$ contains a special subgroup of $G$.
\end{lemma}

\begin{proof}
    Since $d$ is filling, there exists some $rR' \in \spl(d)$ such that
   	\[
   		\bigcap_{\lambda \colon \overline{\reach}_{d}(\{rR'\}) \to R'} \langle d^{t}|_{\underline{s \lambda(s)}} \mid sR' \in \overline{\reach}_{d}(\{rR'\}) \rangle = R
   	\]
   	for all $t \in D$. By \cref{lem:forks in core}, we find some $t$ such that $rR' \in \pi^{\mathrm{ab}}_{R}(\core(M^{t}) \cap R)$. Evidently, also $d^{t} \in \core(M^{t})$ holds.
	By \cref{lem:reach in core}, we may replace $rR'$ by its reach. Now since $\pi^{\mathrm{ab}}_{R}(\core(M^{t}) \cap R)$ is a group, we find $\langle \reach(rR') \rangle \leq \pi^{\mathrm{ab}}_{R}(\core(M^{t} ) \cap R)$, and, applying \cref{lem:reach in core} multiple times, also $\overline{\reach}(rR') \leq \pi^{\mathrm{ab}}_{R}(\core(M^{t} ) \cap R)$.
	
	Let $S$ be a lift of $\overline{\reach}(rR')$ under $\pi^{\mathrm{ab}}_{R}$. For every $s \in S$, we find $d^{t}|_{\underline{s}} = (d^{t})^s|_{\underline{e}} \in \core(M^{t})$, using that $d^t$ and $s$ are contained in $\core(M^t)$. By the choice of $rR'$, the set of these sections generates $R$. Hence $R \leq \core(M^{t})$.
	
	It remains to see that $\langle R \cup \{d^{t}\} \rangle$ is special, i.e.\ that $\langle d^{t}|_{\underline{r}} \mid r \in R\smallsetminus\{e\} \rangle = R$. Since $d$ is filling, this is necessarily true.
\end{proof}


\subsection{Proof of \cref{thm:main positive}} 
\label{sub:proof_of_cref_thm_main_positive}

Finally, we have all the ingredients to {prove the following criterion for inclusion into $\mathcal{MF}$.}

\begin{restatable}{thm}{mainpositive}\label{thm:main positive}
	Let $G = \langle R \cup D \rangle$ be a constant spinal group such that $\mathfrak{SN}(G) = D$. Assume that $D$ is compatible and nilpotent. If there exists some prime power $p^n$ and some surjective homomorphism ${\eta} \colon D \to \Z/p^n\Z$ such that the pre-image of $(\Z/p^n\Z)^\times$ is filling, then $G \in \mathcal{MF}$.
\end{restatable}

\begin{proof}
	Let $M$ be a prodense subgroup of $G$. By \cref{lem:mf no proper prodense} and \cref{thm:projections of proper pd}, it is enough to show that there exists some $u \in X^\ast$ such that the projection $M_u$ is equal to $G_u = G$.
	
	By \cref{prop:prodense conts directed elm}, there exists some $u \in X^\ast$ such that $M_u$ contains an element $d \in {\eta}^{-1}((\Z/p^n\Z)^\times))$. By assumption, $d$ is compatible and filling. Using \cref{lem:special in prodense} we find some element $t \in D$ such that $\{ d^{t} \} \cup R \subseteq \core((M_u)^{t})$. Since $M_u$ is a prodense subgroup of $G$ by \cref{thm:projections of proper pd}, also $(M_u)^{t}$ is a prodense subgroup of $G$. Thus by \cref{prop:all cosets of dirs} we find
	\[
		\pi^{\mathrm{ab}}_{R}(\core((M_u)^{t}) \cap D) = D/D'.
	\]
	Since $D$ is a finite nilpotent group, every lift of $D/D'$ to $D$ generates $D$, hence $\core((M_u)^{t})$ contains $R \cup D$. Thus there is some vertex $w = u {\underline{e}}^n$ such that $((M_u)^{t})_{{\underline{e}}^n} = (M_{w})^{t} = G$, hence $M_w = G$. This concludes our proof.
\end{proof}

Following along the proof above, we see that we need the condition that $D$ is nilpotent only in the very last step. Using the same proof but omitting said step, we obtain the following corollary.

\begin{corollary}\label{cor:special subgroups}
	Let $G$ be a constant spinal group with the same properties as assume in \cref{thm:main positive}, except that $D$ is not necessarily nilpotent. Then every prodense subgroup of $G$ admits a projection containing a special subgroup.
\end{corollary}

Note that the maximal subgroups of infinite index of non-periodic {\v Suni\'k} groups introduced in \cite{FG18}, which are spinal but not constant spinal groups, are not {`special'} in the natural extension of the concept to spinal groups, which can be deduced from the analysis done on the orbital Schreier graphs carried out loc. cit. 



\section{Applications of \cref{thm:main positive}} 
\label{sub:applications_of_cref_thm_main_positive}

{In this section, we prove \cref{thm:primary multi-GGS}, \cref{thm:non mn} and \cref{thm:mf not just-insol}, which we restate at the appropriate point for the convenience of the reader. Furthermore, we provide some additional examples.}

%
%

{As a first application of~\cref{thm:main positive}}, we recover the fact that every periodic prime multi-GGS group~$G$ is contained in $\mathcal{MF}$. In fact, we may {extend the result to} cover a large collection of periodic primary multi-GGS groups.

\primarymultiGGS*

Note that every periodic prime multi-GGS group {vacuously} fulfils the additional condition and is thus covered by the theorem.

\begin{proof}
	Write $H = \langle R \cup \langle d \rangle \rangle$ for the primary GGS group generated by $R$ and $d$ only. This is a special subgroup of $G$. Using \cite[Lemma~2.4]{DFG22}, we may assume that
	\[
		d|_{\underline{r}} = r,
	\]
	otherwise we pass to an $\Aut(X^\ast)$-conjugate GGS group.
	
	Since $G$ is periodic and, in particular, $H$ is periodic, we find that
	\[
		\prod_{r \in R\smallsetminus\{e\}} d|_{\underline{r}} = \id
	\]
	and $\mathfrak{SN}(H) = D$, {using that $r$ generates $R$ and employing the description of periodic GGS groups given in} \cite{Pet23} and \cite{Vov00}. Thus there exists some $y \in X \smallsetminus\{e, r\}$ such that $d|_{\underline{y}} = {r}^m$ also generates $R$; here we use that $R$ is a cyclic group of prime power order. Write $y = {r}^n$ for some $n \in \Z/p^n\Z$.
	
	Now ${r}$ is a forking point for $d$, since the identities $d|_{\underline{r}} = x$, $d|_{\underline{r^n}} = d|_{\underline{y}} = r^{m}$ and $d|_{\underline{(r^m)^{m^{-1}}}} = d|_{\underline{r}} = r$ show that there is a directed path of length {$1$} labelled $(1)$ and a directed path of length $2$ labelled $(n, m^{-1})$ from ${r}$ to itself. Since ${r}$ already generates $R$, {the element} $d$ is filling. Since $R$ is abelian, the directed group $D$ is compatible. Finally, $d$ is of order $p^n$ and thus a filling pre-image of $(\Z/p^n\Z)^\times$ under $\Dab$. Thus by \cref{thm:main positive}, the group $G$ is contained in $\mathcal{MF}$.
\end{proof}

In the next theorem, we address a question/remark raised in \cite{FG18} about the possible maximal subgroups of branch groups{; namely, if every maximal subgroup of a branch group in $\mathcal{MF}$ is necessarily normal. Groups where every maximal subgroup is normal are of importance with regards to the problem of invariable generation, cf.\ \cite{CT23} for a recent study with a view towards branch groups.} We remark that also the construction used to prove \cref{thm:mf not just-insol} provides a proof for \cref{thm:non mn}, but the definition employed in the following proof is easier to deal with and illustrates how to construct varying examples of constant spinal groups inside $\mathcal{MF}$ using \cref{thm:main positive}.

\begin{restatable}{thm}{nonmn}\label{thm:non mn}
	There exists a finitely generated branch group that is in $\mathcal{MF}$ but not in the class $\mathcal{MN}$ of groups with all maximal subgroups normal.
\end{restatable}

\begin{figure}
	\centering
	\begin{tikzpicture}[font=\footnotesize, scale=0.8]
		\node (e_0)	[circle, fill=black, inner sep=0.05cm, label=below:$   R'$] at (-8,0) {};
		\node (e_1)	[circle, fill=black, inner sep=0.05cm, label=below:$t  R'$] at (-6,0) {};
		\node (e_2)	[circle, fill=black, inner sep=0.05cm, label=below:$t^2R'$] at (-4,0) {};

		\path[->] (e_1) edge [loop above]				node			{1}	()
						edge 							node 	[below] {2}	(e_2)
			      (e_2) edge [loop above]				node			{1}	()
						edge [in=45, out=135]			node	[above]	{2}	(e_1);
				  
		\node at (-2,0) {for $d$,};
		
		\node (f_0)	[circle, fill=black, inner sep=0.05cm, label=below:$   R'$] at (0,0) {};
		\node (f_1)	[circle, fill=black, inner sep=0.05cm, label=below:$t  R'$] at (2,0) {};
		\node (f_2)	[circle, fill=black, inner sep=0.05cm, label=below:$t^2R'$] at (4,0) {};
		
		\path[->] (f_1) edge [loop above]				node			{2}	()
						edge 							node 	[below] {1}	(f_2)
			      (f_2) edge [loop above]				node			{2}	()
						edge [in=45, out=135]			node	[above]	{1}	(f_1);
		
		\node at (6,0) {for $d^2$};
	\end{tikzpicture}
	\caption{The graph $\mathcal{R}^{\mathrm{ab}}(d^j)$ as described in the proof of \cref{thm:non mn}.}\label{fig:21 example graph}
\end{figure}
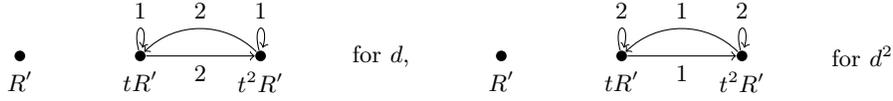

\begin{proof}
	Let $R$ be the semi-direct product $\mathrm{C}_{3} \ltimes \mathrm{C}_7$ given by the presentation
	\[
		\langle a, t \mid a^{7} = t^3 = 1, a^t = a^4 \rangle,
	\]
	i.e.\ the unique non-abelian group of order~$21$. Recall that $R' = \langle [a,t] \rangle^{\mathrm{C}_{3} \ltimes \mathrm{C}_7} = \mathrm{C}_{7}$ and that $R$ is generated by the elements $t^a$ and $t^2$. We define a directed element by
	\[
		d = ({\underline{e}}: d,\, {\underline{R' \smallsetminus\{e\}}}: \id,\, {\underline{tR}}': t^a, \, {\underline{t^2R'}}: t^2).
	\]
	
	It is easy to see that $d$ is compatible, of order $3$, and that $\langle R \cup \langle d \rangle \rangle$ is a constant spinal group. The graphs $\mathcal{R}^{\mathrm{ab}}(d^i)$ for $i \in \{1, 2\}$ are depicted in \cref{fig:21 example graph}; evidently the permutations induced by the $1$-labeled part of the graph are induced by the non-trivial automorphisms of $\mathrm{C}_3$, i.e.\ by the action of $\F_3^\times$ on $\F_3$. It becomes apparent that $tR'$ is a forking point for all $d^i$, using (for example) the following paths
	\begin{align*}
		tR' \xrightarrow{1} tR' \quad&\text{ and }\quad tR' \xrightarrow{2} t^2R' \xrightarrow{2} tR' &\text{ for } d,\\
		tR' \xrightarrow{2} tR' \quad&\text{ and }\quad tR' \xrightarrow{1} t^2R' \xrightarrow{1} tR' &\text{ for } d^2.
	\end{align*}
	Since the reach closure of $tR'$ contains $tR'$, it contains the full commutator factor group. Evidently, every lift $\lambda\colon R/R' \to R'$ provides $\langle d^i|_{\underline{r \lambda(r)}} \mid r \in R/R' \rangle \geq \langle \{(t^a)^i, t^{2i}\} \rangle = R$. Thus $D = \langle d \rangle$ is filling, compatible and nilpotent (indeed, abelian). To apply \cref{thm:main positive}, it remains to check that $\mathfrak{SN}(G) = D$.
	
	{We make use of} \cref{thm:condition for stab nuc}. Thus we have to (partially) compute the map $\Sigma_{d} \colon \mathcal{P}(R) \to \mathcal{P}(R)$. Let $r \in R' = \mathrm{C}_{7}$ be non-trivial. Evidently $d|_{(\underline{r^s)^n}} = \id$ for all $s \in R$ and $n \in {[6]}$. Thus $\Sigma_{d}(r) = \{e\}$. If $r \in tR'$, we find that $r$ has order $3$. Thus $d|_{\underline{r^s}}d|_{\underline{(r^s)^2}} = t^a t^2 \in R'$ for all $s \in R$, similarly $r \in t^2R'$ produces a permutation of the product $d|_{\underline{r^s}} d|_{\underline{(r^s)^2}}$. Since the latter part of the product defining $\Sigma_{\{d\}}(\{r\})$ is always contained in $R'$, we find
	\[
		\Sigma_{d}(r) \subseteq R'
	\]
	for all $r \in R$. Now $\Sigma_{d}(R') = \cup_{r \in R'}\Sigma_{d}(r) = \{e\}$, hence $\Sigma_{d}$ is eventually trivial. By \cref{thm:condition for stab nuc}, $\mathfrak{SN}(G) = D$, and by \cref{thm:main positive}, we find $G \in \mathcal{MF}$.
	
	We now establish that~$G$ is indeed a branch group. To do so, we shall prove that~$G$ is regular branch over the following subgroup
	\[
		B = \langle \{[d, a], [d, t, a], [d, t, t], [d, t, d] \} \rangle^G.
	\]
	This subgroup is of finite index. To see this, consider an element $g = r_0 \prod_{i = 1}^n (d^{\epsilon_i})^{r_i}$ in minimal syllable form, i.e.\ with $n \in \N, \epsilon_i \in \{1, 2\}$ and $r_i \in R$ for $i \in [0, n]$. Every $r_i$ is of the form $a^{\kappa_i}t^{\delta_i}$ for $\kappa_i \in [0,6]$ and $\delta_i \in \{0, 1\}$. Computing modulo $B$, we see that~$a$ commutes with $d$, hence we may assume $\kappa_i = 0$. Thus, two neighbouring syllables are of the form
	\[
		(d^i)^{t^k} (d^j)^{t^l} \equiv_B (d^{i+j})^{t^{k}} [d, t]^{j+l-k}
	\]
	for some $i, j, k, l \in \{1, 2\}$. Since $[d, t]$ is central modulo $B$, we may write
	\[
		g \equiv_B r_0 [d, t]^{\epsilon_0} \prod_{i = 1}^m (d^{\epsilon_i})^{t^{\delta_i}}
	\]
	for some $\epsilon_i$ and $\delta_i \in [0, 2]$, and, more importantly, some $m \leq \lceil n/2 \rceil$, since every two neighbouring syllables produce at most one syllable modulo $B$ using the reduction described above. Continuing that process, we see that $g \equiv_B r_0 [d, t]^{\epsilon_0} (d^{\epsilon_1})^{t^{\delta_1}}$. Since there are only finitely many equivalence classes left and $[d,t]$ has finite order as $\mathfrak{SN}(G) = D$ and \cref{thm:condition for stab nuc}, $B$ is a subgroup of finite index.
	
	To show that $G$ is regular branch over $B$, it is now sufficient to show that for every $s \in \{[d, a], [d, t, a], [d, t, t], [d, t, d]\}$ we find an element $g_s$ in $B$ such that $g_s|_{\underline{e}} = s$ and $g_s|_{\underline{r}} = \id$ for all $r \in R\smallsetminus\{e\}$.
	
	Consider the commutator of $d$ and $(d^2)^{t^2a^4}$. Since both elements stabilise the first {layer}, we may compute their commutator sections-wise, obtaining
	\[
		g = [d, (d^2)^{t^2a^4}] = ({\underline{e}}: [d, t], \, {\underline{ta^{-2}}}: [t^a, d^2], \, {\underline{t^2R'}}: a),
	\]
	using that $[t^2,(t^a)^2] = a$. Now consider
	\[
		[g, d^{t^2}] = ({\underline{e}}: [d, t, t^2], \, {\underline{t^2R'}}: [a, t^a])
	\]
	We take another commutator {to} get rid of the section at ${\underline{e}}$. Computing section-wise as above, we find
	\[
		c = [g, d^{t^2}, d^a] = ({\underline{t^2R'}}: [a, t^a, t^2]).
	\]
	Note that $[a, t^a, t^2] = a^3$. Since
	\[
		g = [d, t^2a^4] [d, a^3t]^{t^2a^4d^2} \equiv_B [d, t^2] [d, t]^{t^2} = 1,
	\]
	the elements $g$ and $c$ are contained in $B$. Now consider the element
	\[
		h = gc^2 = ({\underline{e}}: [d, t], \, {\underline{ta^{-2}}}: [t^a, d^2]) \in B.
	\]
	Now it is a straight-forward computation that
	\begin{align*}
		[h, (d^2)^{t^2}] = ({\underline{e}}: [d,t,t]),\quad
		[h, (c^5)^{t^2}] = ({\underline{e}}: [d, t, a]),\quad\text{and}\quad
		[d, (c^5)^{t^2}] = ({\underline{e}}: [d,a]).
	\end{align*}
	Thus the subgroup $L = \langle \{[d, a], [d,t,t], [d,t,a]\} \rangle^G \leq B$ is geometrically {contained in}~$B$. Consider the element
	\[
		[h, d] = ({\underline{e}}: [d, t, d], \, {\underline{ta^{-2}}}: [t^a, d^2, t^a]) \in B.
	\]
	Since $L$ is geometrically contained in $B$, it is enough to show that $[t^a, d^2, t^a] \in L$ to obtain $(e: [d, t, d]) \in B$, which shows that $B$ is geometrically contained in itself. But
	\[
		[t^a, d^2, t^a] \equiv_L [t, d^2, t^a] \equiv_L [t, d^2, a^4][t, d^2, t]^{a^4} \in L.
	\]
	Thus $G$ is a branch group.
	
	It remains to observe that $H = \langle \{t\} \cup \St_G(1) \rangle$ is a maximal subgroup of $G$ and that $t^a \notin H$. Thus $G \notin \mathcal{MN}$.
\end{proof}

It is possible to prove an even stronger statement. All previous examples of branch groups inside $\mathcal{MF}$ are \emph{just-insoluble}, i.e.\ insoluble groups such that every proper quotient is soluble, while on the contrary, almost all known examples of branch groups outside $\mathcal{MF}$ do not have this property. We complete the picture by constructing not just-insoluble branch groups inside $\mathcal{MF}$.

\mfnotjustinsol*

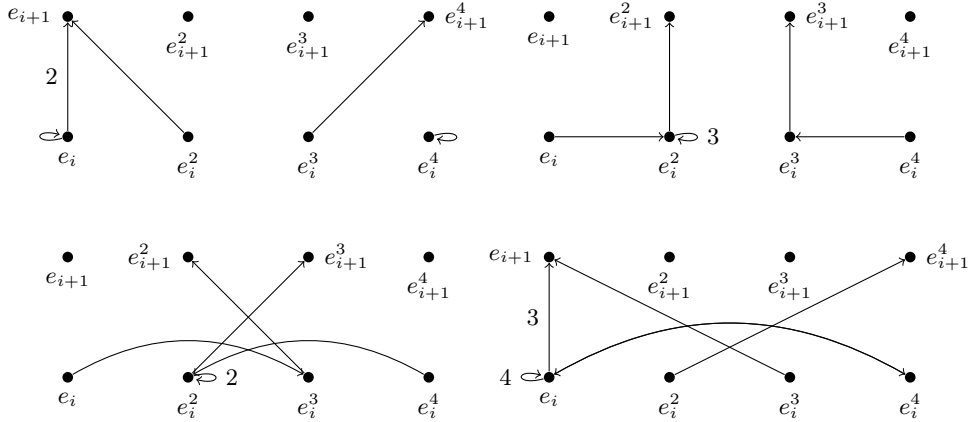
\begin{figure}
	\centering
	\begin{tikzpicture}[font=\footnotesize, scale=0.8]
		\node (e_1)	[circle, fill=black, inner sep=0.05cm, label=below:$e_i  $] at (0,2) {};
		\node (e_2)	[circle, fill=black, inner sep=0.05cm, label=below:$e_i^2$] at (2,2) {};
		\node (e_3)	[circle, fill=black, inner sep=0.05cm, label=below:$e_i^3$] at (4,2) {};
		\node (e_4)	[circle, fill=black, inner sep=0.05cm, label=below:$e_i^4$] at (6,2) {};

		\node (g_1)	[circle, fill=black, inner sep=0.05cm, label=left: $e_{i+1}  $] at (0,4) {};
		\node (g_2)	[circle, fill=black, inner sep=0.05cm, label=below:$e_{i+1}^2$] at (2,4) {};
		\node (g_3)	[circle, fill=black, inner sep=0.05cm, label=below:$e_{i+1}^3$] at (4,4) {};
		\node (g_4)	[circle, fill=black, inner sep=0.05cm, label=right:$e_{i+1}^4$] at (6,4) {};

		\path[->] (e_1) edge [loop left]				node	[left]	{}	()
		edge 							node 	[left] 	{2}	(g_1)
		(e_2) edge 							node	[below]	{}	(g_1)
		(e_3)	edge 							node	[below]	{}	(g_4)
		(e_4)	edge [loop right]				node	[right]	{}	();
		
		\begin{scope}[shift={(8,0)}]
		\node (e_1)	[circle, fill=black, inner sep=0.05cm, label=below:$e_i  $] at (0,2) {};
		\node (e_2)	[circle, fill=black, inner sep=0.05cm, label=below:$e_i^2$] at (2,2) {};
		\node (e_3)	[circle, fill=black, inner sep=0.05cm, label=below:$e_i^3$] at (4,2) {};
		\node (e_4)	[circle, fill=black, inner sep=0.05cm, label=below:$e_i^4$] at (6,2) {};
		
		\node (g_1)	[circle, fill=black, inner sep=0.05cm, label=below:$e_{i+1}$]   at (0,4) {};
		\node (g_2)	[circle, fill=black, inner sep=0.05cm, label=left: $e_{i+1}^2$] at (2,4) {};
		\node (g_3)	[circle, fill=black, inner sep=0.05cm, label=right:$e_{i+1}^3$] at (4,4) {};
		\node (g_4)	[circle, fill=black, inner sep=0.05cm, label=below:$e_{i+1}^4$] at (6,4) {};

		\path[->] (e_1)	edge 							node 	[left]  {}	(e_2)
			      (e_2) edge [loop right]				node	[right]	{3}	()
				  	    edge 							node	[below]	{}	(g_2)
				  (e_3)	edge 							node	[below]	{}	(g_3)
				  (e_4)	edge 							node	[right]	{}	(e_3);
		\end{scope}
		
		\begin{scope}[shift={(0,-4)}]
			\node (e_1)	[circle, fill=black, inner sep=0.05cm, label=below:$e_i  $] at (0,2) {};
			\node (e_2)	[circle, fill=black, inner sep=0.05cm, label=below:$e_i^2$] at (2,2) {};
			\node (e_3)	[circle, fill=black, inner sep=0.05cm, label=below:$e_i^3$] at (4,2) {};
			\node (e_4)	[circle, fill=black, inner sep=0.05cm, label=below:$e_i^4$] at (6,2) {};
		
			\node (g_1)	[circle, fill=black, inner sep=0.05cm, label=below: $e_{i+1}  $] at (0,4) {};
			\node (g_2)	[circle, fill=black, inner sep=0.05cm, label=left:$e_{i+1}^2$] at (2,4) {};
			\node (g_3)	[circle, fill=black, inner sep=0.05cm, label=right:$e_{i+1}^3$] at (4,4) {};
			\node (g_4)	[circle, fill=black, inner sep=0.05cm, label=below:$e_{i+1}^4$] at (6,4) {};

			\path[->] (e_1)	edge [bend left]				node 	[left]  {}	(e_3)
				      (e_2) edge [loop right]				node	[right]	{2}	()
					  	    edge 							node	[below]	{}	(g_3)
					  (e_3)	edge 							node	[below]	{}	(g_2)
					  (e_4)	edge [bend right]				node	[right]	{}	(e_2);
		\end{scope}
		\begin{scope}[shift={(8,-4)}]
			\node (e_1)	[circle, fill=black, inner sep=0.05cm, label=below:$e_i  $] at (0,2) {};
			\node (e_2)	[circle, fill=black, inner sep=0.05cm, label=below:$e_i^2$] at (2,2) {};
			\node (e_3)	[circle, fill=black, inner sep=0.05cm, label=below:$e_i^3$] at (4,2) {};
			\node (e_4)	[circle, fill=black, inner sep=0.05cm, label=below:$e_i^4$] at (6,2) {};
		
			\node (g_1)	[circle, fill=black, inner sep=0.05cm, label=left: $e_{i+1}  $] at (0,4) {};
			\node (g_2)	[circle, fill=black, inner sep=0.05cm, label=below:$e_{i+1}^2$] at (2,4) {};
			\node (g_3)	[circle, fill=black, inner sep=0.05cm, label=below:$e_{i+1}^3$] at (4,4) {};
			\node (g_4)	[circle, fill=black, inner sep=0.05cm, label=right:$e_{i+1}^4$] at (6,4) {};

			\path[->] (e_1) edge [loop left]				node	[left]	{4}	()
							edge 							node 	[left] 	{3}	(g_1)
							edge [bend left]				node			{}	(e_4)
				      (e_2) edge 							node	[below]	{}	(g_4)
					  (e_3)	edge 							node	[below]	{}	(g_1)
					  (e_4)	edge [bend right]				node	[right]	{}	(e_1);
		\end{scope}	
		
	\end{tikzpicture}
	\caption{Parts of the graph $\mathcal{R}^{\mathrm{ab}}(d^j)$ for $j = 1$ to $4$ top left to bottom right as described in the proof of \cref{thm:mf not just-insol}. Not all edges and vertices are drawn. Unlabelled edges are thought to be labelled with~$1$.}\label{fig:rab for thm mf not just-insol}
\end{figure}

\begin{proof}
	Let $n \in \N_{>4}$ and let $m(n)$ be the minimal size of a generating set for $\Alt(n)$ consisting only of $5$-cycles, e.g.\ $m(5) = 2$; such a generating set always exists, since the normal subset of all $5$-cycles must generate a normal subgroup. Let $S = \{s_1, \dots, s_{m(n)}\}$ be such a set. Define $R = \mathrm{C}_5^{m(n)} \times \Alt(n)$ and set $E = \{e_1, \dots, e_{m(n)}\}$ to be a generating set for $\mathrm{C}_5^{m(n)}$. Evidently $R' = \Alt(n)$ and $R/R' = \mathrm{C}_5^{m(n)}$.
	
	Define a directed automorphism $d$ by
	\[
		d|_{\underline{r}} = \begin{cases}
			s_i e_i		&\text{ if } r \in e_iR' \text{ for some } i \in [1, m(n)],\\
			e_{i+1}		&\text{ if } r \in e_i^2R' \text{ for some } i \in [1, m(n)],\\
			e_{i+1}^4	&\text{ if } r \in e_i^3R' \text{ for some } i \in [1, m(n)],\\
			e_{i}^4		&\text{ if } r \in e_i^4R' \text{ for some } i \in [1, m(n)],\\
			d			&\text{ if } r = e,\\
			\id			&\text{ otherwise,}
		\end{cases}
	\]
	where one reads indices modulo $m(n)$, i.e.\ $e_{m(n)+1} = e_1$. Set $D = \langle d \rangle$ and consider the constant spinal group $G = \langle R \cup D \rangle$. Evidently $D$ is compatible and a cyclic group of order $5$. A diagram sketching $\mathcal{R}^{\mathrm{ab}}(d^j)$ for $j \in {[4]}$ can be found in \cref{fig:rab for thm mf not just-insol}.
	
	From it (or directly from the description of $d$ above) we may easily read of that $e_iR'$ is a forking point of $d$ and $d^4$ for all $i \in {[m(n)]}$, and that $e_i^2R'$ is a forking point of $d^2$ and $d^3$ for all $i \in {[m(n)]}$, using the paths
	\begin{align*}
		e_iR' \xrightarrow{1} e_iR' \quad&\text{ and }\quad e_iR' \xrightarrow{2} e_{i+1}R' \xrightarrow{2} \dots \xrightarrow{2} e_{i-1}R' \xrightarrow{2} e_iR' &\text{ for } d,\\
		e_i^2R' \xrightarrow{3} e_i^2R' \quad&\text{ and }\quad e_i^2R' \xrightarrow{1} e_{i+1}^2R' \xrightarrow{1} \dots \xrightarrow{1} e_{i-1}^2R' \xrightarrow{1} e_iR' &\text{ for } d^2,\\
		e_i^2R' \xrightarrow{2} e_i^2R' \quad&\text{ and }\quad e_i^2R' \xrightarrow{4} e_{i+1}^2R' \xrightarrow{4} \dots \xrightarrow{4} e_{i-1}^2R' \xrightarrow{4} e_iR' &\text{ for } d^3,\\
		e_iR' \xrightarrow{4} e_iR' \quad&\text{ and }\quad e_iR' \xrightarrow{3} e_{i+1}R' \xrightarrow{3} \dots \xrightarrow{3} e_{i-1}R' \xrightarrow{3} e_iR' &\text{ for } d^4.
	\end{align*}
	Observing the seconds path described above for every directed element $d^j$, it is apparent that $\reach(e_1R')$ and, respectively, $\reach(e_1^2R')$ contain $\{e_iR' \mid i \in [1, m(n)]\}$, a generating set for $R/R'$. Thus $\overline{\reach}(e_1R') = R/R'$ and $\overline{\reach}(e_1^2R') = R/R'$.

	Since $d$ is not only compatible, but the map $\Delta_{d^j}$ is {actually} constant on cosets of $R'$, there is only one lift to consider, which yields $\{(s_ie_i)^j \mid i \in [1, m(n)]\} \cup \{e_i^{\pm j} \mid i \in [1, m(n)]\}$, which is a generating set for $R$ for all $j \in [1,4]$. Thus every non-trivial power of $d$ is filling.
	
	To apply \cref{thm:main positive}, we have yet to show that $\mathfrak{SN}(G) = D$. We aim to use \cref{thm:condition for stab nuc} and (partially) compute $\Sigma_{d} \colon \mathcal{P}(R) \to \mathcal{P}(R)$. Since the sections of $d$ only depend on the $R'$-coset of the vertex, we find $\Sigma_{d}(r) = \Sigma_{d}(s)$ for all $r \equiv_{R'}s$. The projection of $\{s^n \mid s \in r^R, n \in [1, \ord(r)) \}$ under $\pi^{\mathrm{ab}}_{R}$ is a cyclic subgroup with the trivial element removed. It is easily derived from the definition of $d$ that a product running (perhaps multiple times) over vertices lying above any cyclic subgroup of $R/R'$ is equal to a power of some $s_i$ modulo $R'$. Furthermore, the individual factors generate a subgroup of the form $\langle s_i \rangle \times \langle e_i \rangle \times \langle e_{i+1} \rangle$, which is abelian. Thus the product is indeed equal to a power of $s_i$. Hence
	\[
		\Sigma_{d}(r) = \begin{cases}
			\langle s_i \rangle &\text{ if }r \in \langle e_iR' \rangle \smallsetminus\{R'\} \text{ for some } i \in [1, m(n)]\\
			\{e\}&\text{ otherwise}.
		\end{cases}
	\]
	In particular, $\Sigma^2_{d}(r) = \{\id\}$ for all $r \in R$. Thus $\mathfrak{SN}(G) = D$ by \cref{thm:condition for stab nuc}, and by \cref{thm:main positive} we find $G \in \mathcal{MF}$. It remains to notice that $G/\St_G(1) \cong R$ has the insoluble quotient $R/C_5^{m(n)}\cong \Alt(n)$.
	
	It remains to see that at least one of the groups is branch; we show that $G$ is branch over the subgroup $B = \langle [D, R] \rangle^G$ in case $n = 5$. The procedure is similar as in the last proof. Consider the element
	\[
		s = [d, d^{e_1^4e_2}] = ({\underline{e_1R'}}: [s_1, s_2]) \in B.
	\]
	Since the element $[s_1, s_2]$ normally generates $\Alt(5)$ and conjugation by $d$ and $d^{e_1^4e_2}$ correspond to the simultaneous conjugation of the block of sections equal to $[s_1, s_2]$ by $s_1$ and $s_2$, respectively, we find that $\overline{r} = ({\underline{e_1R'}}: r) \in B$ for all $r \in \Alt(5)$. Conjugating by appropriate elements of $R$, we find elements of the form $(R': r)$ and $({\underline{e_2R'}}: r) \in B$. Thus $({\underline{e}}: [d, t]) \in B$ for all $t \in \Alt(5)$, i.e. $\langle [D, \Alt(5)]\rangle^G$ is geometrically contained in $B$. Furthermore, by multiplication of $d$ with elements of the forms described above we find
	\begin{align*}
		\tilde{d} = ({\underline{e}}: d,\, &{\underline{e_1R'}}: e_1,\, {\underline{e_1^2R'}}: e_2,\, {\underline{e_1^3R'}}: e_2^4,\, {\underline{e_1^4R'}}: e_1^4,\\
		&{\underline{e_2R'}}: e_2,\, {\underline{e_2^2R'}}: e_1,\, {\underline{e_2^3R'}}: e_1^4,\, {\underline{e_2^4R'}}: e_2^4\,) \in dB.
	\end{align*}
	One computes
	\[
		[\tilde{d}, \tilde{d}^{e_1}] = ({\underline{e}}: [d, e_1], \, {\underline{e_1^4}}: [e_1^4, d]) \in B
	\]
	and furthermore
	\begin{align*}
		[\tilde{d}, \tilde{d}^{e_1}, \tilde{d}^{e_1s_1}] &= ({\underline{e}}: [d, e_1, e_1]) \in B,\\
		[\tilde{d}, \tilde{d}^{e_1}, \tilde{d}^{e_2s_1}] &= ({\underline{e}}: [d, e_1, e_2]) \in B,
	\end{align*}
	hence $\langle [d, e_1, e_1] \rangle^G$ is contained geometrically in $B$. {Furthermore}, we find 
	\[
		b = {[\tilde{d}, \tilde{d}^{e_1}, \tilde{d}^{s_1}]^{-1} =} ({\underline{e_1^4}}: [e_1^4, d, e_1^4]^{-1}) \in B,
	\]
	{hence}
	\[
		[\tilde{d}, \tilde{d}^{e_1}, \tilde{d}]b = ({\underline{e}}: [d, e_1, d]) \in B.
	\]
	Similarly, one obtains $[d, e_2, e_2], [d, e_2, e_1]$ and $[d, e_2, d]$. In view of \cref{lem:passage to normal closure}, we find
		\[
		\gamma_3(G) = \langle \Alt(5) \cup \{[d, r] \mid r \in \Alt(5)\} \cup \{[d, e_i, e_j], [d, e_i, d] \mid i, j \in \{1, 2\}\} \rangle^G
		\]
		geometrically contained in $B$. Thus
		\[
		[\tilde{d}, \tilde{d}^{e_1}] \equiv_B ({\underline{e}}: h, \, {\underline{e_1^4}}: h) = b \in B
		\]
		where we put $h = [d, e_1]$. It remains to compute
		\[
		b(b^{-1})^{e_1}b^{e_1^2}(b^{-1})^{e_1^3}b^{e_1^4} \equiv_B ({\underline{e}}: h^2) \in B,
		\]
		and, proceeding similarly, $({\underline{e}}: [d, e_2]) \in B$. Thus $B$ is geometrically contained in itself. The quotient $G/B$ is an image of $R \times D$, a finite group, so $B$ is also of finite index and $G$ is a branch group.
\end{proof}

The construction in the previous proof is readily seen to be generalised by replacing $\Alt(n)$ by any group $K$ that can be generated by elements of prime power order for a prime $p > 3$; in case $K$ is not perfect, set the sections of the element $d$ to be constant on the cosets of $K$ instead of constant only on the cosets of the commutator subgroup. Thus, a plethora of finite groups may be realised as a quotient of a constant spinal group in $\mathcal{MF}$. On the other hand, \cref{thm:main positive} is not applicable to constant spinal groups with perfect rooted groups. Thus the following question seems natural:

{\begin{question}
	Does there exist a finitely generated branch group $G$ in $\mathcal{MF}$ such that the group $G|^u$ is perfect for every $u \in \Aut(X^\ast)$?
\end{question}}

As a final {application of \cref{thm:main positive}}, we introduce a constant spinal group with non-abelian directed group that is contained in $\mathcal{MF}$, contrasting the previous examples and demonstrating the full potential of \cref{thm:main positive}.

\begin{example}[A constant spinal group with non-abelian directed group contained in $\mathcal{MF}$]
	Denote by $R$ the Heisenberg group over the field with $5$ elements, generated by the two elements $a$ and $b$. We define a directed generator --- that is similar to the one defined in the previous example --- by the $1$-labelled part of the graph $\mathcal{R}^{\mathrm{ab}}(d_1)$, given by the following diagram,
	
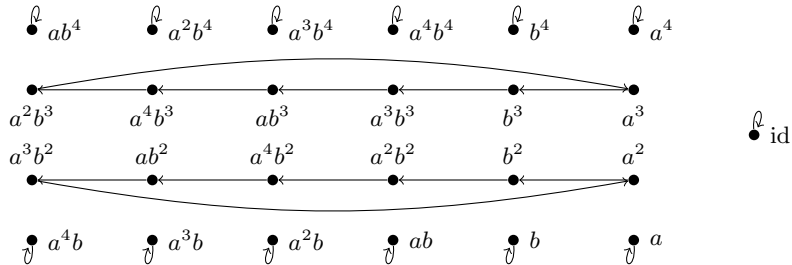
\begin{center}
		\begin{tikzpicture}[font=\footnotesize, scale=0.8, rotate=90]
			\node (a_1)	[circle, fill=black, inner sep=0.05cm, label=right:$a  $] at (0,0)   {};
			\node (a_2)	[circle, fill=black, inner sep=0.05cm, label=above:$a^2$] at (1,0)   {};
			\node (a_3)	[circle, fill=black, inner sep=0.05cm, label=below:$a^3$] at (2.5,0) {};
			\node (a_4)	[circle, fill=black, inner sep=0.05cm, label=right:$a^4$] at (3.5,0) {};
		
			\node (b_1)	[circle, fill=black, inner sep=0.05cm, label=right:$b  $] at (0,2)   {};
			\node (b_2)	[circle, fill=black, inner sep=0.05cm, label=above:$b^2$] at (1,2)   {};
			\node (b_3)	[circle, fill=black, inner sep=0.05cm, label=below:$b^3$] at (2.5,2) {};
			\node (b_4)	[circle, fill=black, inner sep=0.05cm, label=right:$b^4$] at (3.5,2) {};
		
			\node (ab_1)	[circle, fill=black, inner sep=0.05cm, label=right:$ab    $] at (0,4)   {};
			\node (ab_2)	[circle, fill=black, inner sep=0.05cm, label=above:$a^2b^2$] at (1,4)   {};
			\node (ab_3)	[circle, fill=black, inner sep=0.05cm, label=below:$a^3b^3$] at (2.5,4) {};
			\node (ab_4)	[circle, fill=black, inner sep=0.05cm, label=right:$a^4b^4$] at (3.5,4) {};
		
			\node (a2b_1)	[circle, fill=black, inner sep=0.05cm, label=right:$a^2b  $] at (0,6)   {};
			\node (a2b_2)	[circle, fill=black, inner sep=0.05cm, label=above:$a^4b^2$] at (1,6)   {};
			\node (a2b_3)	[circle, fill=black, inner sep=0.05cm, label=below:$a  b^3$] at (2.5,6) {};
			\node (a2b_4)	[circle, fill=black, inner sep=0.05cm, label=right:$a^3b^4$] at (3.5,6) {};
		
			\node (a3b_1)	[circle, fill=black, inner sep=0.05cm, label=right:$a^3b  $] at (0,8)   {};
			\node (a3b_2)	[circle, fill=black, inner sep=0.05cm, label=above:$a  b^2$] at (1,8)   {};
			\node (a3b_3)	[circle, fill=black, inner sep=0.05cm, label=below:$a^4b^3$] at (2.5,8) {};
			\node (a3b_4)	[circle, fill=black, inner sep=0.05cm, label=right:$a^2b^4$] at (3.5,8) {};
		
			\node (a4b_1)	[circle, fill=black, inner sep=0.05cm, label=right:$a^4b  $] at (0,10)   {};
			\node (a4b_2)	[circle, fill=black, inner sep=0.05cm, label=above:$a^3b^2$] at (1,10)   {};
			\node (a4b_3)	[circle, fill=black, inner sep=0.05cm, label=below:$a^2b^3$] at (2.5,10) {};
			\node (a4b_4)	[circle, fill=black, inner sep=0.05cm, label=right:$a  b^4$] at (3.5,10) {};
		
			\node (tr)		[circle, fill=black, inner sep=0.05cm, label=right:$\id$] at (1.75,-2) {};
		
			\path[->]	(tr)	edge	[loop right]	();
		
			\path[->]	(a_1)	edge	[loop left]		()
						(a_4)	edge	[loop right]	()
						(a_2)	edge					(b_2)
						(a_3)	edge					(b_3);
		
			\path[->]	(b_1)	edge	[loop left]		()
						(b_4)	edge	[loop right]	()
						(b_2)	edge					(ab_2)
						(b_3)	edge					(ab_3);
					
			\path[->]	(ab_1)	edge	[loop left]		()
						(ab_4)	edge	[loop right]	()
						(ab_2)	edge					(a2b_2)
						(ab_3)	edge					(a2b_3);
					
			\path[->]	(a2b_1)	edge	[loop left]		()
						(a2b_4)	edge	[loop right]	()
						(a2b_2)	edge					(a3b_2)
						(a2b_3)	edge					(a3b_3);
					
			\path[->]	(a3b_1)	edge	[loop left]		()
						(a3b_4)	edge	[loop right]	()
						(a3b_2)	edge					(a4b_2)
						(a3b_3)	edge					(a4b_3);
		
			\path[->]	(a4b_1)	edge	[loop left]		()
						(a4b_4)	edge	[loop right]	()
						(a4b_2)	edge	[bend right=10]	(a_2)
						(a4b_3)	edge	[bend left=10]		(a_3);
		\end{tikzpicture}
	\end{center}
	where we have dropped the $R'$ behind each representative of a coset, and by stating that if $d_1|_{\underline{x}} \equiv_{R'} a^i b^j$, then $d_1|_{\underline{x}} = a^ib^j$. In fact, any lift such that the product over the sections in $R'\smallsetminus\{e\}$ vanishes would do. Note that the columns of the graph above correspond to the cyclic subgroups of $R/R'$. Additionally we define another directed automorphism $d_2$ satisfying $d_2|_{\underline{x}} = \id$ for all $x \in R\smallsetminus(\{e\} \cup aR' \cup a^4R')$, and $d_2|_{\underline{as}} = b$ and $d_2|_{\underline{a^4s}} = b^4$ for all $s \in R'$. We view this directed automorphism as some error possibly perturbing $d_1$. We now consider the constant spinal group $G = \langle R \cup D \rangle$ with $D = \langle \{d_1,d_2\} \rangle$. Note that $D$ is isomorphic to the Heisenberg group over $\F_5$.
	
	We now check {the requirements of~\cref{thm:main positive}} as usual. We may consider the map $D\colon \Z/5\Z$ with kernel $\langle d_2 \rangle^D = \langle d_2 \rangle D'$, hence we have to show that the pre-image $\bigcup_{k = 1}^4 d_1^k\langle d_2 \rangle D'$ of the set of generators is filling. Crucially, we do not need to show that $d_2$ is filling (indeed, $d_2$ is not filling). We provide the following paths, proving that its first vertices $b^k$ are forking points for $d_1^k$, with $k \in {[4]}$:
	\begin{align*}
		b	\xrightarrow{1} b	&\quad\text{and}\quad	b	\xrightarrow{2} a^2b^2
															&\xrightarrow{1} a^4b^2
															\xrightarrow{1} ab^2
															\xrightarrow{1} a^3b^2
															\xrightarrow{1} a^2
															\xrightarrow{1} b^2
															\xrightarrow{3} b		&\quad\text{for }d_1,\\
		b^2	\xrightarrow{3} b^2	&\quad\text{and}\quad	b^2	\xrightarrow{1} a^4b^4
															&\xrightarrow{1} a^3b^3
															\xrightarrow{1} a^2b
															\xrightarrow{1} a^4b^2
															\xrightarrow{1} a^2b^4
															\xrightarrow{1} a^4b^3
															\xrightarrow{1} a^4b\\
															&&\xrightarrow{1} a^3b^2
															\xrightarrow{1} a^4
															\xrightarrow{2} b^2
															\xrightarrow{1} b		&\quad\text{for }d_1^2,\\
		b^3	\xrightarrow{2} b^3	&\quad\text{and}\quad	b^3	\xrightarrow{1} a^4b^4
															&\xrightarrow{1} a^2b^2
															\xrightarrow{1} a^2b
															\xrightarrow{1} ab^3
															\xrightarrow{1} a^2b^4
															\xrightarrow{1} ab^2
															\xrightarrow{1} a^4b\\
															&&\xrightarrow{1} a^2b^3
															\xrightarrow{1} a^4
															\xrightarrow{3} b
															\xrightarrow{1} b^3		&\quad\text{for }d_1^3,\\
		b^4	\xrightarrow{1} b^4	&\quad\text{and}\quad	b^4	\xrightarrow{3} a^3b^3
															&\xrightarrow{1} a^4b^2
															\xrightarrow{1} a^4b^3
															\xrightarrow{1} a^3b^2
															\xrightarrow{1} a^3
															\xrightarrow{1} b^2
															\xrightarrow{2} b^4		&\quad\text{for }d_1^4.
	\end{align*}
	These elements passed through these paths are already generating $R/R'$, and since $R$ is nilpotent, every lift of a generating set of $R/R'$ generates the full group $R$. Furthermore, notice that the paths for $d_1$ and its inverse do not pass through the set $\{a, a^4\}$ at all, while the paths for $d_1^2$ and its inverse exit that set over an edge labeled~$2$ or~$3$. Consider that the graphs $\mathcal{R}^{\mathrm{ab}}(d_1^j)$ and $\mathcal{R}^{\mathrm{ab}}(d_1^jd_2^k)$ for $j \in {[4]}$ and $k \in [0, 4]$ differ only in the edges labeled $\pm i \bmod 5$ passing out of $a^i$, since $(d_1^jd_2^k)|_{\underline{x}} = d_1^j|_{\underline{x}}$ for all $x \notin \{a, a^4\}$. Thus the paths as sketched above also witness that $b^j$ is a forking point for all directed elements $d_1^jd_2^k$.
	
	It remains to show that $\mathfrak{SN}(G) = D$, for which we compute $\Sigma_{\{d_1, d_2\}}$. The columns of the diagram sketched above correspond to the (non-trivial elements of the) cyclic subgroups of $R/R'$. Looking at the products over the columns, it is easy to see that they are contained in $R'$. This yields $\Sigma_{\{d_1, d_2\}}(r) \subseteq R'$ for all $r \in R'$. Since the sections of $d_1$ and $d_2$ at elements within $R'\smallsetminus\{\id\}$ are trivial, we have $\Sigma_{\{d_1, d_2\}}(R') = \{\id\}$. Thus $\mathfrak{SN}(G) = D$ and by \cref{thm:main positive} we find $G \in \mathcal{MF}$.
	
	Note that one may easily replace the Heisenberg group over $\F_5$ by its equivalent over $\F_p$ for all primes $p > 3$ by filling up the graphs for $d_1$ and $d_2$ with arrows pointing to $\id$.
\end{example}


\section{Maximal subgroups of infinite index in branch groups} 
\label{sec:maximal_subgroups_of_infinite_index_in_branch_groups}

In this section we establish \cref{thm:maximal subgroups in layered constant spinal subgroups} {and furthermore employ it} to explicitly construct maximal subgroups of infinite index in certain branch groups, closing some gaps in {the} understanding of the class of branch groups \emph{outside} $\mathcal{MF}$.

To begin with, we establish some preliminary results on the structure of {constant spinal groups and their} special subgroups, which are interesting in their own right. An automorphism $g \in \Aut(T)$ is called \emph{finitary} if there exists a finite set $Y \subset T$ of vertices such that the label $g|^u$ is trivial for all $u \in T\smallsetminus Y$, or, equivalently, such that $g|_u = \id$ for all but finitely many $u \in T$. Every rooted automorphism is evidently finitary.

\begin{proposition}[{Classification of finitary elements}]\label{lem:labels of finitary elements}
	Let $G$ be a constant spinal group. Let $g \in G$ be a finitary automorphism and $u \in X^* \smallsetminus\{\varnothing\}$. Then $g|^u \in R'$.
\end{proposition}

\begin{proof}
	We make the following claim: If $v \in X^\ast$ is such that that $g|_{vx} \equiv_{G'} s_{vx} \in R$ for all $x \in X$, then $g|^{vx} \in R'$ for all $x \in X$ and $g|_v \equiv_{G'} s_v \in R$.
	
	Assume the claim holds. Let ${v} \in X^\ast$ be the {predecessor} vertex of $u$, which exists since $u$ is not the root. Since $g$ is finitary, there exist vertices $v_1, \dots, v_k$ for some $k \in \N$ such that $g|_{v_jx} \in R$ for all $j \in [1, k]$ and $x \in X$ (i.e.\ fulfilling the assumption of the claim) and such that every infinite ray {passing through} $v$ contains one vertex $v_j$. By the claim, we may replace each $v_j$ by its {predecessor} vertex. Since it also fulfils the assumption of the claim, {we may continue to pass to predecessors} until we find that $g|_{vx} \equiv_{G'} s_{vx} \in R$ for all $x \in X$ and hence $g|^{u} \in R'$.
	
	It remains to proof the claim. Let $v$ be such that $g|_{vx} \equiv_{G'} s_{vx} \in R$ for all $x \in X$. Write $g|_v = r_0 \prod_{i = 1}^n d_i^{r_i}$ in syllable form. By our assumption, for every $x \in X$ we find
	\[
		\prod_{i \in J_x(g|_v)} d_i \equiv_{D'} \id.
	\]
	It follows that $g|_v \equiv_{G'} g|^v \prod_{x \in X} \prod_{i \in J_x(g|_v)} d_i = g|^v \in R$. Thus we have established the second part of the claim.
	
	{One checks that given $d \in D'$ and $x \in X\smallsetminus\{\underline{e}\}$, we have $d|_x \in R'$. Using this fact, we compute the following, where we put $K = \langle \St_G(1)\cup R' \rangle$,
	\begin{align*}
		(g|_v)|_x &\equiv_{K} \prod_{y \in X\smallsetminus\{x\}} \prod_{i \in J_y(g|_v)} d_i|_{y^{-1}x}\\
		&= \prod_{y \in X\smallsetminus\{x\}} \left(\prod_{i \in J_y(g|_v)} d_i\right)|_{y^{-1}x} \in R'.
	\end{align*}
	Thus} we find that $g|^{vx} = g|_{vx} \bmod \St_G(1)$ is contained in $R'$.
\end{proof}

Note that, clearly, $g|^\varnothing$ may be any element of $R$. {Furthermore, the existence of constant spinal groups that are regular branch over their commutator subgroup and possess non-abelian rooted groups (cf.\ \cref{eg:groups_without_csp}) shows that finitary elements with arbitrary portrait adhering to the statement of the previous lemma may exist.} It is also worth to mention that the lemma above (and the following lemma, building on it) do not hold in case of constant spinal groups with rooted groups acting non-regularly: For example, if $d$ is a directed element and $r$ an element fixing the distinguished point $e \in X$ but fulfilling $r.x = y$ for some $x \neq y \in X$. Then $[d, r] = d^{-1} d^r$ fulfils $[d, r]|_{\underline{e}} = d^{-1}d = \id$. All other sections are rooted, hence $[d,r]$ is finitary, but $[d, r]|_y = d|_x^{-1} d|_y$ may, depending on $d$, take any value in $R$.

\begin{proposition}[Layer-climbing]\label{lem:layer-climbing lemma}
	Let $G = \langle R \cup D \rangle$ be a constant spinal group, and let $H = \langle R \cup T \rangle$ be a special subgroup of $G$. Assume that $H$ is regular branch over its commutator subgroup $H'$. Let $n \in \N_0$ be a non-negative integer. If an element $g \in G$ satisfies that all its $n$\textsuperscript{th} layer sections are contained in $H$, then $g \in H$.
\end{proposition}

\begin{proof}
	The {statement} is trivially true for $n = 0$. Otherwise, we may reduce to the case $n = 1$ {in the following way}. If the lemma is true for $n = 1$, and $g \in G$ fulfils $g|_{X^n} {\subseteq} H$ for some $n > 1$, we may consider all sections at vertices in $uX$ for some $u \in X^{n-1}$. Since they are all contained in $H$,  also $g|_u$ is contained in $H$, and the case $n$ reduces to the case $n - 1$. Thus assume that $g|_x \in H$ for all $x \in X$.
	
	Since $H$ is a constant spinal group, its factor commutator group is isomorphic to $R/R' \times T/T'$, therefore $g|_x \equiv_{H'} r_x d_x$ for some $r_x \in R$ and $d_x \in T$. Since $H'$ is geometrically contained in $H$, there exists some $h_1 \in H$ such that $(gh_1)|_x = r_x d_x$. Now consider the element
	\[
		h_2 = \prod_{x \in X} (d_x)^{x^{-1}},
	\]
	where the order of the product may be chosen arbitrary. Evidently $h_2|_x \equiv_{H'} s_xd_x$ for some $s_x \in R$. Thus $gh_1h_2^{-1}$ fulfils $(gh_1h_2^{-1})|_x \in R$ for all $x \in X$. By \cref{lem:labels of finitary elements}, we find $(gh_1h_2^{-1})|_x = r_xs_x^{-1} \in R'$. But $R' \leq H'$ is geometrically contained in $H$, hence there exists $h_3 \in H$ such that $h_3|_x = s_xr_x^{-1}$. It follows that the product $gh_1h_2^{-1}h_3$ is rooted. Since $R \leq H$, we find $g \in H$.
\end{proof}

From here on, it becomes feasible to {describe} the maximal subgroups among special prodense subgroups {and prove \cref{thm:maximal subgroups in layered constant spinal subgroups}, we we restate for the convenience of the reader. Note that while this class of subgroups is clearly greatly restricted, this can be used to show that maximal subgroups of infinite index exist.}

\maximalsubgroupsinlayeredconstantspinalsubgroups*

\begin{proof}
	Since maximality is preserved under conjugation, we may assume $H = \langle R \cup D \rangle$.
	
	If the subgroup $T$ is not maximal in $D$, any (automatically special) subgroup $\tilde{T}$ of $G$ generated by $R$ and $T < \tilde{T} < D$ clearly contains $H$. Notice that a special subgroup $K = \langle R \cup L \rangle$ is proper if and only if the corresponding subgroup $L$ of $D$ is a proper subgroup of $D$. This follows directly from the description of the nucleus of a constant spinal group, which fulfils $\mathfrak{N}(K) \equiv_R L$, which is equal to $\mathfrak{N}(G) \equiv_R D$ if and only if $D = L$. Thus $\langle R \cup \tilde{T} \rangle$ is a proper supergroup of $H$.
	
	Now assume that $T$ is a maximal subgroup of $D$. Let $g \in G \smallsetminus H$ and set $K = \langle H \cup \{g\} \rangle$, which is itself a prodense subgroup of $G$. We aim to show that $K$ is in fact equal to $G$. By the description of the nucleus of $G$ we find some $n \in \N_0$ such that $g|_{u} \equiv_R d_u \in D$ for all $u \in X^n$. By \cref{lem:layer-climbing lemma}, we may assume that there exists at least one $u \in X^n$ such that $d_u \notin T$, otherwise we would find $g \in H$. Since $H$ is spherically transitive, self-similar and fractal, there exists some $h \in H$ such that $gh \in \st_K(u)$ and $(gh)|_u = g|_u = d_u$. Thus $K_u = \langle H \cup \{d_u\} \rangle$, using the fact that $R \leq H_u = H$ since $H$ is a constant spinal group. But clearly $\langle H \cup \{d_u\} \rangle = \langle R \cup T \cup \{d_u\} \rangle = G$, since $T$ is a maximal subgroup of $D$. Now \cref{thm:projections of proper pd} implies that $K = G$, hence $H$ is a maximal subgroup of $G$.
\end{proof}

It remains to actually construct prodense special subgroups. In the case of a cyclic group of prime order as the rooted group, this is impossible: {here}, the constant spinal groups are precisely the multi-GGS groups. {These groups} do not permit prodense special (necessarily also multi-GGS) subgroups, cf.\ \cite{KT18}, {which can also be seen as} a special case of our next proposition.

\begin{proposition}\label{prop:necessary condition for density}
	Let $G = \langle R \cup D \rangle$ be a constant spinal group. Let $H = \langle R \cup T \rangle$ be a special subgroup. If $H$ is prodense, then $T^D (R^G \cap D) = D$. In particular, if $R$ or~$D$ is nilpotent, there are no proper prodense special subgroups.
\end{proposition}

\begin{proof}
	If $H$ is prodense, we have that $H^G = H^GH = G$. Furthermore, $H^G = T^DR^G$, so by Dedekind's law
	\[
		D = D \cap G = D \cap T^D R^G = (D \cap R^G) T^D.
	\]
	In view of \cref{thm:factor commutator groups of constant spinal groups}, we find $R^G \cap D \leq D'$. If $R$ is nilpotent, so is $D$, which we now assume. Since the derived subgroup of a nilpotent group is supplemented only by the full group, $D = T^D (R^G \cap D)\leq T^DD'$ show that $T^D = D$. But every maximal subgroups of a nilpotent group is normal, hence $T = D$ and $H = G$.
\end{proof}

Unfortunately, this rules out our method as a source of counterexamples for Passman's second conjecture on $p$-groups. On the other hand, it is not difficult to produce prodense special subgroups in layered constant spinal groups. It is even more straight-forward in the case that $D$ is perfect; which is a consequence of the following result.

\begin{proposition}\label{prop:perfectness and the csp}
	Let $G = \langle R \cup D \rangle$ be a layered constant spinal group. Then $G$ has the congruence subgroup property if and only if $G$ is perfect.
\end{proposition}

\begin{proof}
	The commutator subgroup is a proper normal subgroup of finite index. Assume that $G$ has the congruence subgroup property. Thus exists some $n \in \N$ such that $\St_G(n) \leq G'$. Since $G$ is layered, we have $\St_G(n) = G \times \dots \times G$. Thus there exists some element $g_d \in \St_G(n)$ such that $g_d|_u = \id$ for all $u \in X^n \smallsetminus\{ v \}$ and $g_d|_v = d$ for some $v \in X^n$ and all $d \in D$. If there exists some $d \notin D'$, we have ${\Dab}(g_d) = dD' \neq D'$, using the map ${\Dab}$ defined in the proof of \cref{thm:factor commutator groups of constant spinal groups}. Thus $g_d \notin G'$ which contradicts our assumption $\St_G(n) \leq G'$ therefore implies $D = D'$. At the same time, $R$ must be perfect; otherwise, the infinitely iterated wreath product of a non-perfect group would be topologically finitely generated, which is impossible, cf.\ \cite{Bon10}. Thus $G/G'$ is trivial and $G$ is perfect. 
	
	On the other hand, if $G$ is perfect, using the fact that it is layered, we find that
	\[
		\Rist_G(n)' = (G \times \dots \times G)' = G \times \dots \times G = \St_G(n)
	\]
	for all $n \in \N$. Since every normal subgroup of any branch group necessarily contains the commutator subgroup of some rigid stabiliser, this implies that $G$ has the congruence subgroup property.
\end{proof}

Note that this allows us to easily construct large classes of groups without the congruence subgroup property, as demonstrated in the following example. {Finding natural examples of such groups is not a wholly trivial task, however, following the first example given by Pervova, cf.\ \cite{Per07}, now many more are known; see \cite{BSZ12, GS23, Ski20} for other examples.}

\begin{example}[Groups without the congruence subgroup property]\label{eg:groups_without_csp}
	Let $R$ be a finite non-abelian simple group. Let $a$ and $b$ be two distinct and non-trivial elements of $R$ such that the equations $(a^{-1}b)^2 \neq 1$, $a^2 \neq b$ and $b^2 \neq a$ hold. Let $S = \{s_0, s_1\}$ be a generating set for $R$. Define a directed automorphism by
	\[
		d = ({\underline{e}}: d, \, {\underline{a}}: s_0, \, {\underline{b}}: s_1).
	\]
	Since $S$ generates $R$, the group $G = \langle R \cup \langle d \rangle \rangle$ is a constant spinal group. It is evidently not perfect, since $\langle d \rangle$ is cyclic. Thus by \cref{prop:perfectness and the csp} it is enough to show that $G$ is layered. This can be seen by the following computation,
	\[
		[d, d^{ba^{-1}}] = ({\underline{e}}: [d, d|_{ba^{-1}}], \, {\underline{a}}: [s_0, s_1], \, {\underline{b}}: [s_1, d|_{ba^{-1}b}]).
	\]
	By choice of the elements $a$ and $b$, the product $ba^{-1}$ is neither trivial nor equal to $a$ or $b$, the same holds for $ba^{-1}b$. Thus $[d, d^{ba^{-1}}] = ({\underline{a}}: [s_0, s_1])$. Since $[s_0, s_1]$ normally generates the simple non-abelian group $R$, we find $R$ geometrically {contained in $G$, cf.\ \cref{lem:passage to normal closure}.} Since $d \equiv_{R\times \dots \times R} ({\underline{e}}: d)$, the whole group $G$ is geometrically {contained in itself}, i.e.\ $G$ is layered.
\end{example}

We now come back to the construction of prodense special subgroups. We shall use the following idea for our construction: Let $G = \langle R \cup T \rangle$ be a layered constant spinal group. {By \cref{lem:layered_csp_justinfinite}, the group $G$ is just-infinite and its profinite and normal topologies coincide.} Let $D$ be a perfect supergroup of $T$ in the group of directed automorphisms. Evidently $\langle R \cup D \rangle$ is also a layered and a constant spinal group. Since $D$ is perfect, the group has the congruence subgroup property and since the (special) subgroup $H = \langle R \cup T \rangle$ is also layered, it is congruence-dense, thus profinitely dense, and thus prodense.

To construct many such pairs, we {provide} two preliminary lemmata.

\begin{lemma}\label{Lemma: generators R one}
	Let $R$ be a finitely generated group and $S$ a minimal generating set for $R$. If $rs_1$ and $rs_2$ belong to $S$ for some non-trivial $r \in R$ and $s_1, s_2 \in S$ with $s_1 \neq s_2$, then $rs_1 = s_2$, $rs_2 = s_1$, and $r^2 = e$.
\end{lemma}

\begin{proof}
	Let us denote $t_1 = rs_1$ and $t_2 = rs_2$. Then 
	\begin{align*}
		r = s_1^{-1}t_1 = s_2^{-1}t_2.
	\end{align*}
	Since $S$ is a minimal generating set, we deduce that $t_1, t_2 \in \{s_1,s_2\}$. By hypothesis, $r \neq 1$, hence we find $s_2 = t_1 = rs_1$ and $s_1 = t_2 = rs_2$. This yields $r^2s_1 = s_1$, i.e.\ $r$ has order 2. 
\end{proof}

\begin{lemma}\label{Lemma: generators R two}
	Let $R$ be a finitely generated group. Then there exists a minimal generating set $S$ for $R$ with the following property: For any $s_1,s_2 \in S$, if $s_2s_1^{-1}s_2 = s_1$, then $[s_1,s_2] = e$.
\end{lemma}

\begin{proof}
	Suppose that both $s_1$ and $s_2$ are involutions. Then $s_2s_1^{-1}s_2 = s_1$ implies
	\begin{align*}
		e = s_1^{-1}s_2s_1^{-1}s_2 = s_1^{-1}s_2^{-1}s_1s_2 = [s_1,s_2].
	\end{align*}
	Otherwise, we may assume without loss of generality that $s_2$ is not an involution. By \cref{Lemma: generators R one}, the equality $s_2s_1^{-1}s_2 = s_1$ implies that $s_2s_1^{-1}$ is an involution. Furthermore, $(S \smallsetminus \{s_2\}) \cup \{s_2s_1^{-1}\}$ is again a minimal generating set for $R$. Following this procedure and using that $S$ is finite, we can ensure that there exists a minimal generating set $\widetilde{S}$ for $R$ with the following property: For any $s_1$ and $s_2$ in $\widetilde{S}$, if $s_2s_1^{-1}s_2 = s_1$, then both $s_1$ and $s_2$ are involutions, which by the first part of this proof implies that $[s_1,s_2]=1$.
\end{proof}

\begin{example}[A family of branch groups {outside $\mathcal{MF}$}]\label{eg:class_of_examples}
	Let $R$ be a finite perfect group. We define the following family of directed elements: For every $s, k \in R\smallsetminus\{e\}$ consider the directed element
	\begin{align*}
		d_{s,k} = ({\underline{e}}: d_{s,k},\, {\underline{s}}:k).
	\end{align*}
	Then $D = \langle \{ d_{s,k} \mid s,k \in R\setminus\{e\} \}\rangle \cong R \times\stackrel{|R|-1}{\cdots}\times R$ is a perfect group. Since a group generated by the union of two perfect groups is perfect, the group $G = \langle R \cup D \rangle$ is a perfect constant spinal group.
	
	Now, using \cref{Lemma: generators R one} and \cref{Lemma: generators R two}, there exists a minimal generating set $S$ for $R$ and a set $\Omega = \{(s_1,s_2) \mid s_1,s_2 \in S\}$ with the following properties:
	\begin{enumerate}
		\item $\langle \{ [s_1,s_2] \mid (s_1,s_2) \in \Omega \} \rangle^R = R$;
		\item for every pair $(s_1,s_2) \in \Omega$, $s_2s_1^{-1}s \notin S$ for all $s \neq s_1$ in $S$.
	\end{enumerate}
	Consider the directed element
	\begin{align*}
		d = ({\underline{e}}:d,\, {\underline{s}}: s \, \mid s \in S) \in D.
	\end{align*}
	Fix $(s_1,s_2) \in \Omega$, and put $r = s_2s_1^{-1}$. Then 
	\begin{align*}
		d^{r^{-1}} =({\underline{r}}: d,\, {\underline{rs}}: s \, \mid s\in S)
	\end{align*}
	and 
	\begin{align*}
		[d^{r^{-1}},d]= ({\underline{s_2}}: [d|_{s_1}, d|_{s_2}] = [s_1, s_2]).
	\end{align*}
	We deduce that every special subgroup of $G$ containing $d$ is conjugate to a layered group, and, since $G$ has the congruence subgroup property by \cref{prop:perfectness and the csp} {and is just infinite by \cref{lem:layered_csp_justinfinite}}, is furthermore prodense in $G$. Thus let $T$ be a maximal subgroup of $D$ containing $d$. By \cref{thm:maximal subgroups in layered constant spinal subgroups}, we deduce that $H = \langle R \cup T \rangle$ is a special proper prodense subgroup of $G$, thus $H$ is a maximal subgroup of $G$ of infinite index.
\end{example}

Note that these {groups} are all acting on trees with rather large valency, e.g.\ starting with the $60$-regular tree $\Alt(5)^\ast$. The example of {a non-$\mathcal{MF}$ branch group by} Francoeur and Garrido, on the other hand, acts of the binary tree. To the authors' knowledge, { it is unknown if} there exist{s} a branch group outside of $\mathcal{MF}$ acting on the $3$-regular tree. Note that the methods of Bondarenko \cite{Bon10} do not apply to the $3$-regular tree.

\arbitrarylongchainsofmaximalsubgroups*

\begin{proof}
	Let $R = \operatorname{Alt}(n)$ be an alternating group of degree $n \geq 7$ and odd. Fix the $7$-cycle $s = (1\,2\,3\,4\,5\,6\,7)$ and the double $3$-cycle $t = (1\,2\,3)(4\,5\,6)$. Consider the following directed elements
	\begin{align*}
		a &= ({\underline{e}}: a, \, {\underline{t}}: (1\, 2\, 3)),\\
		b &= ({\underline{e}}: b, \, {\underline{t}}: (3\, 4\, \cdots\, n)),\\
		c &= ({\underline{e}}: c, \, {\underline{s}}: (1\, 2\, 3)),\\
		d_{k} &= \begin{cases}
			({\underline{e}}: d_{k}, \, {\underline{s}}:(3\ 4\ \cdots\ {k})) \qquad &\text{if $7\leq {k}\leq n$ is odd,}\\
			({\underline{e}}: d_{k}, \, {\underline{s}}:(3\ 4\ \cdots\ {k}-1)) \qquad &\text{if $7<{k}< n$ is even},
			\end{cases}\\
		\widetilde{d_{k}} &= \begin{cases}
			\id \qquad &\text{if $7 \leq n \leq n$ is odd,}\\
			({\underline{e}}: \widetilde{d_{k}}, \, {\underline{s}}: (4\ \cdots\ {k})) \qquad &\text{if $7 < {k} < n$ is even}.
		\end{cases}
	\end{align*}
	Note that  $\langle \{a,b,c,d_{k}, \widetilde{d_k}\}\rangle \cong \operatorname{Alt}(n) \times\operatorname{Alt}(k)$ is a perfect group, for all $k \in [7, n]$. Therefore
	\begin{align*}
		G_k = \langle R \cup \{a,b,c,d_{k}, \widetilde{d_k}\} \rangle
	\end{align*}
	is also perfect for every $k \in [7,n]$. Furthermore, the group $G_k$ is a layered constant spinal group: Indeed, if we consider the element $r = ts^{-1}$, and compute 
	\begin{align*}
		[(cb)^r, cb] = (s: [(cb)|_t, (cb)|_s]) = (s: [(34\cdots n), (123)])).
	\end{align*}
	We easily deduce that $R$ is geometrically contained in $G_k$. This implies that we may find all elements of the form $(e:a)$, $(e:b)$, $(e:c)$,  $(e: d_k)$, and $(e: \widetilde{d_k})$ in $G_k$ for all odd $m \in [7, k]$, showing that $G_k$ is layered. Using {the fact} that $G_k$ is perfect and layered, we find
	\begin{align*}
		G_7 \leq_{\text{prodense}} G_8 \leq_{\text{prodense}} \cdots \leq_{\text{prodense}} G_n.
	\end{align*}
	Since the directed group of $G_{k-1}$ and $G_k$ are isomorphic to $\operatorname{Alt}(k-1) \times \operatorname{Alt}(n)$ and $\operatorname{Alt}(k) \times \operatorname{Alt}(n)$, respectively, the inclusion $G_{k-{1}} < G_k$ is proper, and using \cref{thm:maximal subgroups in layered constant spinal subgroups} we see that  $G_{k-1}<_{\text{max}} G_{k}$, and {the} chain above is indeed a chain of stepwise maximal subgroups of infinite index.
\end{proof}

{As another application, we show that there exists a periodic branch group outside $\mathcal{MF}$, in contrast to the example of a non-$\mathcal{MF}$ branch group of Garrido and Francoeur. Note that the fact that said groups is non-periodic is an integral part of their approach to show that it contains maximal subgroups of infinite index.}

\torsionbranchgroupoutsidemf*

\begin{proof}
	We construct an infinite family of examples. Let $n \geq 5$ be an integer fulfilling. Put $R = \Alt(n+2)$. Let $q_0$ and $q_1$ be two elements in $R$ with the following properties,
	\begin{enumerate}
		\item $q_i^R \cap \Alt(n+1) = \varnothing$ for $i \in \{0, 1\}$,
		\item $\ord(q_i) > 2$ for $i \in \{0, 1\}$,
		\item $q_0$ and $q_1$ are not conjugate,
		\item $\langle q_0 \rangle \cap \langle q_1 \rangle = \{e\}$,
		\item $(q_0q_1^{-1})^2 \neq e$,
		\item $q_0^2[q_0, q_1] \neq e$,
		\item $q_1^{q_0}\neq q_1^{\pm1}$.
	\end{enumerate}
	For example, if $n$ is even, such a pair is given by
	\[
		q_0 = (1\,2\, \dots\, n)((n+1) \, (n+2)) \quad\text{and}\quad q_1 = (1\,2\, \dots\, (n-1))(n \, (n+1) \, (n+2)).
	\]
	For odd $n$, one may choose
	\begin{align*}
		q_0 &= (1\,2\, \dots\, (n-3))((n-2)\,(n-1)\,n)((n+1) \, (n+2)) \quad\text{and}\\ q_1 &= (1\,2\, \dots\, (n-2))((n-1) \,n) ((n+1) \, (n+2)).
	\end{align*}
	{However, there are many more possible choices of pairs fulfilling the above properties.}
	
	Let $\sigma \in R$ be such that $\langle \Alt(n)^\sigma \cup \Alt(n) \rangle = R$, for example the double transposition~$(1\, (n+1))(2\, (n+2))$.\\
	Using property (3), we define a directed element $\varphi(r)$ for every $r \in \Alt(k)$ with $k\in \N$
	\begin{align*}
		\varphi(r) &= ({\underline{e}}: \varphi(r),\, {\underline{q_0}}: r,\, {\underline{q_0^{-1}}}: r^{-1},\, {\underline{q_1}}: r^\sigma,\, {\underline{q_1^{-1}}}: (r^\sigma)^{-1}).
	\end{align*}
	For any $\Alt(k)$ with $k\in \N$, it is apparent that $\varphi$ is an isomorphism between $\Alt(k)$ and $\varphi(\Alt(k))= \{\varphi(r) \mid r \in \Alt(k)\}$. From properties (2),(4), (5), (6) and (7) we obtain that
	\[
		[\varphi(r), \varphi(r)^{q_0q_1^{-1}}] = ({\underline{q_1}}: [r^\sigma,r])
	\]
	for every $r \in \Alt(k)$. Set $D = \varphi(\Alt(n+1))$ and $T = \varphi(\Alt(n))$. Note that $H = \langle R \cup T \rangle$ is a special subgroup of $G = \langle R \cup D \rangle$. Furthermore, since $\sigma$ does not centralise $\Alt(n)$, we find some $r \in \Alt(n)$ such that $[\varphi(r), \varphi(r)^{q_0q_1^{-1}}]$ is not trivial. Using \cref{lem:passage to normal closure}, we find that $R$ is geometrically contained in $H$. Thus $H$ and $G$ are layered and since $\Alt(n)$ is a maximal subgroup of $\Alt(n+1)$, \cref{thm:maximal subgroups in layered constant spinal subgroups} implies that $H$ is a maximal subgroup of $G$.	
	It remains to show that $G$ is periodic, for which we employ \cref{thm:condition for stab nuc}. Let $r \in R$ be contained in a conjugate of $\Alt(n+1)$. Since $r^R \cap \{q_0, q_0^{-1}, q_1, q_1^{-1}\} = \varnothing$ by property (1), we see that $d|_{s} = \id$ for all $s \in \langle r \rangle^R \smallsetminus\{e\}$ and all $d \in D$. Thus $\Sigma_{D}(r) = \{ e \}$. On the other hand, if $r \in R$ is any element, by property (2) $\langle r \rangle^R$ contains at most one of the sets $\{q_0, q_0^{-1}\}$ and $\{q_1, q_1^{-1}\}$. If it contains neither, it has empty intersection with both. If it contains $\{q_i, q_i^{-1}\}$ for either $i=0$ or $i=1$, then
	\[
		\Sigma_D(r) = \langle \{ d|_{\underline{q_i}}d|_{\underline{q_i^{-1}}} \mid d \in D\} \rangle \cdot \langle \{ d|_{\underline{q_i}},d|_{\underline{q_i^{-1}}} \mid d \in D\} \rangle' \leq \Alt(n+1)^{\sigma^i}.
	\]
	But since $\Sigma_D$ maps every $r$ contained in a conjugate of $\Alt(n+1)$ to $e$, we find $\Sigma^2_D(r) = e$ for all $r \in R$. Thus, $G$ is periodic by \cref{thm:condition for stab nuc}.
\end{proof}

Finally, we {construct an example completing} the picture regarding the interplay between the congruence subgroup property and membership in $\mathcal{MF}$ for branch groups. It is known that groups in $\mathcal{MF}$ might or might not exhibit the congruence subgroup property (for example the GGS groups and {some of} the EGS groups, cf. \cite{KT18,Per07}), while previous examples of branch groups outside $\mathcal{MF}$ are all known to possess said property.

\branchgroupoutsidemfwithoutcsp*

\begin{proof}
	Let $R = \Alt(5)$ be the alternating group on five points. {Put} $c_1 = (1 \, 2 \, 3)$ and $c_2 = (3 \, 4 \, 5)$, and furthermore $b_1 = (1 \, 2)(3 \, 4)$ and $b_2 = (1 \, 5)(3 \, 4)$. Define two directed automorphisms by
	\begin{align*}
		t = ({\underline{e}}: t, \, {\underline{c_1}}: c_1, \, {\underline{c_2}}: c_2) \quad\text{ and }\quad
		d = ({\underline{e}}: d, \, {\underline{c_1}}: b_1, \, {\underline{c_2}}: b_2).
	\end{align*}
	Clearly $\langle \{c_i, b_i\} \rangle \cong \Alt(4)$ for $i \in {[2]}$. Beyond that, the {pairs} $(c_1, c_2)$ and $(b_1, b_2)$ generate the diagonal subgroup in $\Alt(4) \times \Alt(4)$, hence $\langle \{d, t\}\rangle = D \cong \Alt(4)$. Set $G = \langle R \cup D \rangle$. This is evidently a constant spinal group. Furthermore, $\langle t \rangle = T \leq D$ gives rise to a special subgroup $H = \langle R \cup T \rangle$.
	
	We now show that $H$ is layered. Consider that $ c_2c_1^2, c_2c_1^2 c_2 \notin \{\id, c_1, c_2\}$, hence
	\[
		[t, t^{c_2c_1^2}] = ({\underline{e}}: [t, t|_{c_2c_1^2}] = \id, \, {\underline{c_1}}: [c_1, d|_{c_2}] = [c_1, c_2], \, {\underline{c_2}}: [c_2, d|_{c_2c_1^2c_2}]).
	\]
	Since $[c_1, c_2]$ normally generates $\Alt(5)$, we find $R$ geometrically {contained in} $H$. Since $t \equiv_{R \times \dots \times R} (e: t)$, also $T$ is geometrically {contained in} $H$; the same holds for~$D$ and $G$. Thus both groups are layered and in particular branch. By \cref{prop:perfectness and the csp}, the group~$G$ does not possess the congruence subgroup property, since it is not perfect. Note that therefore, even though $H$ and $G$ are both layered {(whence profinitely dense subgroups are prodense by \cref{lem:layered_csp_justinfinite})}, it is not necessarily the case that $H$ is prodense. However, we now {demonstrate} that it is indeed prodense.
	
	First observe that $G' H = G$, since $D' T = D$. Every normal subgroup of $G$ contains the derived subgroup of a rigid layer stabiliser $\Rist_G(n)$. Since $G$ is layered, the rigid and usual layer stabilisers coincide, hence every normal subgroup geometrically contains a subgroup of the form $L_n = G' \times \dots \times G'$ with the number of copies of $G'$ {being equal to the cardinality of $X^n$}. Since $H$ is layered, it geometrically contains $\St_H(n) = H \times \dots \times H$ with the same number of copies. Now since $G'H = G$, we find that $L_nH$ contains $G \times \dots \times G = \St_G(n)$, for an appropriate integer $n$. But $\St_G(n)H = G$, since $H$ is congruence-dense in $G$, thus $L_nH = G$ also holds. Since every normal subgroup contains some $L_n$, the group $H$ is prodense. Since $T$ is a maximal subgroup of $D$, it is indeed a maximal subgroup of infinite index in $G$ by \cref{lem:mf no proper prodense}.
\end{proof}
	

\providecommand{\bysame}{\leavevmode\hbox to3em{\hrulefill}\thinspace}
\providecommand{\MR}{\relax\ifhmode\unskip\space\fi MR }
\providecommand{\MRhref}[2]{%
  \href{http://www.ams.org/mathscinet-getitem?mr=#1}{#2}
}
\providecommand{\href}[2]{#2}

\end{document}